\documentclass[nonblindrev]{informs3_no_remark} 

\OneAndAHalfSpacedXI 

\usepackage{natbib}
 \bibpunct[, ]{(}{)}{,}{a}{}{,}%
 \def\bibfont{\small}%

 \usepackage{tikz,calc}
 \usetikzlibrary{decorations.markings}
 \usepackage{soul}
 \setul{}{0.75pt}
 \setuloverlap{1.5pt}

\usepackage{subfigure}
\usepackage{mwe}
\usepackage{lscape,booktabs,longtable}
\usepackage{amsmath}
\usepackage[pdftex,colorlinks=true,urlcolor=blue,citecolor=blue,pdfstartview=FitH]{hyperref}
\usepackage{comment}
\usepackage{enumerate}
\usepackage[ruled,vlined]{algorithm2e}
\usepackage[graphicx]{realboxes}

\newcommand{\E}[1]{{\mathrm{E}\left[#1\right]}}

\newcommand{\PP}[1]{\mathrm{P}\left( #1 \right)}

\newcommand{\edit}{{}}

\TheoremsNumberedThrough     
\ECRepeatTheorems

\EquationsNumberedThrough    

\begin{document}





\RUNAUTHOR{Daw}

\RUNTITLE{Conditional Uniformity and Hawkes Processes}

\TITLE{
Conditional Uniformity and Hawkes Processes
}
\ARTICLEAUTHORS{%
\AUTHOR{Andrew Daw}
\AFF{University of Southern California Marshall School of Business} 

} 



\ABSTRACT{Classic results show that the Hawkes self-exciting point process can be viewed as a collection of temporal clusters, where exogenously generated initial events give rise to endogenously driven descendant events. This perspective provides the distribution of a cluster's size through a natural connection to branching processes, but this is irrespective of time. Insight into the chronology of a Hawkes process cluster has been much more elusive. Here, we employ this cluster perspective and a novel adaptation of the random time change theorem to establish an analog of the conditional uniformity property enjoyed by Poisson processes. Conditional on the number of epochs in a cluster, we show that the transformed times are jointly uniform within a particular convex polytope. Furthermore, we find that this polytope leads to a surprising connection between these continuous state clusters and parking functions, discrete objects central in enumerative combinatorics and closely related to Dyck paths on the lattice. In particular, we show that uniformly random parking functions constitute hidden spines within Hawkes process clusters. This yields a decomposition that is valuable both methodologically and practically, which we demonstrate through application to the popular Markovian Hawkes model and through proposal of a flexible and efficient simulation algorithm. 
}
\KEYWORDS{Hawkes processes, parking functions, conditional uniformity, cluster duration, self-excitement, compensators, time change, Dyck paths} 


\maketitle


%


\section{Introduction}

 The hallmark of the Hawkes self-exciting point process is that each event generates its own stream of offspring events, so that the history of the process becomes an endogenous driver of its own future activity \citep{hawkes1971spectra}. In this way, the model can be viewed as a collection of \emph{clusters} of events, where each cluster is initiated by some exogenous activity and then filled by the progeny descending from that original event. This structure appears throughout the many applications of the Hawkes process, observed in the virality of social media and the contagion of financial risk, and employed in the prediction of seismological activity and the spread of COVID-19 \citep[e.g.][]{ogata1998space,ait2015modeling,bertozzi2020challenges,nickel2020learning}.\footnote{See also: \url{https://ai.facebook.com/blog/using-ai-to-help-health-experts-address-the-covid-19-pandemic/}} This cluster-based perspective was first provided by \citet{hawkes1974cluster}, and it has served as a cornerstone for analysis of the stochastic process. A particularly useful consequence is that, through a connection to Poisson branching processes, we are granted comprehensive time-agnostic insight into the model. That is, the distribution of the size of the cluster, meaning the number of events it contains, is well understood and available in closed form. In this paper, we address a question that is elementary, related, yet evasive: how long will a cluster last?

The duration of a Hawkes process cluster, meaning the time elapsed from the first epoch to the last, must be at least as important in application as the cluster size; examples of longing to know when an endogenously driven activity will cease surely cannot be far from mind. However, despite the ease of access to the distribution of the cluster size, this answer has remained elusive. While formalizing the cluster-based definition, \citet{hawkes1974cluster} bounded the mean cluster duration and obtained an integral equation that must be satisfied by its cumulative distribution function (CDF), but acknowledged that it would only ``be solved in principle by repeated numerical integration.'' \citet{moller2005perfect} tightened the bound on the mean of the duration and showed that the duration itself could be stochastically dominated by an exponential random variable under some tail conditions on the excitation function, but the authors also remarked that the duration distribution is ``unknown even for the simplest examples of Hawkes processes.'' They too develop a numerical approximation of the CDF, and this informs their approach for introducing a perfect simulation procedure for the Hawkes process, which is the primary focus of that work. \citet{chen2021perfect} is also devoted to perfect simulation of the Hawkes process, particularly in steady-state, and applies this to a single server queue driven by Hawkes arrivals. By comparison to \citet{moller2005perfect}, \cite{chen2021perfect} leverages exponential tilting to produce a more efficient algorithm. Because exponential tilting with respect to the cluster duration is challenging, the author instead constructs an upper bound of the duration that is almost surely longer than the duration itself. It is important that we note that simulation of Hawkes process clusters is actually a subroutine of the \citet{chen2021perfect} algorithm \citep[and likewise for][]{moller2005perfect}, and this is explicitly described in the extension to multivariate Hawkes processes in \citet{chen2020perfect}. In fact, all of these works simulate the cluster through the structure provided by the \citet{hawkes1974cluster} definition. 

Several other works have sought headway into the cluster duration outside the context of simulation. \citet{graham2021regenerative} also bounds the duration of the cluster, invoking the tail conditions from \citet{moller2005perfect} to provide inequalities for the mean and for the CDF. Here, the author uses the bounds to prove regenerative properties for the Hawkes process as a whole. Relatedly, \citet{costa2020renewal} also show\edit{s} renewal results, but only for excitation functions with a bounded support condition. Additionally, \citet{reynaud2007some} provide\edit{s} non-asymptotic tail bounds of the cluster duration, and \citet{bremaud2002rate} bound\edit{s} the tail of the duration in order to bound the rate of convergence to equilibrium in the Hawkes process overall. However, none of these prior works characterize the mean duration explicitly, let alone its distribution.

Here, we address the fundamental pursuit of the cluster duration through a surprising connection from Hawkes processes to parking functions, a family of random objects that are rooted in enumerative combinatorics.  In a comprehensive survey of combinatorial results, \citet{yan2015parking} remarked that the parking function lies ``in the center of combinatorics and appear[s] in many discrete and algebraic structures.'' In this paper, we find that parking functions are also the hidden spines of Hawkes process clusters. While this bridge from discrete to continuous space may be unexpected, the parking function itself is truly not a lonely concept. Originally introduced in \citet{konheim1966occupancy}, the classic context is an ordered vector of preferences over a row of $k$ parking spaces such that, if $k$ drivers  proceed left to right and take only their preferred space or the next available to the right of it, all cars will be able to park. In enumerative combinatorics, this simple concept has been connected to many other interesting objects. For example, \citet{riordan1969ballots} linked parking functions and labeled trees, while \citet{stanley1997parking} connected parking functions and non-crossing partitions. \citet{stanley2002polytope} then extended this to plane trees and to the associahedron, as well as to a family of polytopes closely related to the one that arises here. Parking functions have also been of use in the analysis of polynomials, such as in \citet{carlsson2018proof}. 

Highly relevant to our following study is the bijection between parking functions and labeled Dyck paths, or, equivalently, between sorted parking functions and Dyck paths \citep[see e.g.~Section 13.6 in][]{yan2015parking}. Many other combinatorial relationships are available in great detail from \citet{yan2015parking}. Quite recently, the parking function has received considerable attention as a random object; for example by \citet{diaconis2017probabilizing} and~\citet{kenyon2021parking}. 
Motivated by the idea of a uniformly random parking function (from the collection of all those considering the same number of spaces), \citet{diaconis2017probabilizing}  explore\edit{s} the joint distribution of the full vector of preferences and the marginal distribution of the first preference. Additionally, the authors also conduct an asymptotic study of parking functions as the number of cars (or, equivalently, spaces) grows large, yielding a connection to Brownian excursion processes in the limit. In the combinatorics literature, there are many generalizations of the parking function in which the cardinality of spaces exceeds the cardinality of preferences, and \citet{kenyon2021parking} adapt\edit{s} this generalization to the stochastic model. These generalized parking functions are again uniformly random within each combination of space and preference sizes, and \citet{kenyon2021parking} stud\edit{ies} both the marginal distribution of a single coordinate for the generalized parking function and the covariance between two coordinates for the classical notion.


Our connection between Hawkes processes and parking functions joins a family of work descend\edit{ing} from the random time change theorem for point processes. The classical random time change theorem \citep[e.g. Theorem 7.4.I of][]{daley2003introduction} provides that if the integral of the process intensity is evaluated at the arrival epochs of a simple, adapted (not necessarily Poisson) point process, these transformed points will form a unit\edit{-}rate Poisson process. This integral of the intensity is called the compensator. This beautiful result can be traced back to \citet{meyer1971demonstration}, with the corresponding characterization of the Poisson on the half-line originating with \citet{watanabe1964discontinuous}. There is a deep theory surrounding this \edit{transformation}, such as the elegant martingale proofs given in \citet{bremaud1975extension} or \citet{brown1988simple}. A rich collection of these ideas is available in \citet{bremaud1981point}, as are many demonstrations of the potency of martingales for point process models. However, it has since been seen that this result need not require the most advanced of techniques, as \citet{brown2002time} showed that this attractive idea can be explained with elementary arguments, essentially reducing the proof to calculus. We will make use of this simplicity here. 

In addition to \edit{its} elegance, the random time change \edit{theorem} also \edit{holds great value}. For example, in the simulation of point processes, if one can compute and invert the compensator function, Algorithm 7.4.III of \citet{daley2003introduction} shows that one can obtain a point process with the corresponding intensity function by transforming arrival epochs of a unit\edit{-}rate Poisson process according to the inverse of the compensator. For the Hawkes process, the application of this idea can be traced to \citet{ozaki1979maximum}, and inverting the compensator is essentially the idea behind the exact simulation procedure in \citet{dassios2013exact}, even if it is not described as such.  As we will see, our focus on the Hawkes cluster sets us apart. The simulation method that arises from our analysis differs from these Hawkes compensator inversion methods by sampling the points collectively from a polytope, rather than iteratively in the sample path. In statistics, the random time change theorem leads to methods for evaluating the fit of point process models on data \citep[see, e.g.,][]{ogata1988statistical,brown2005statistical,kim2014call}. The idea of these techniques is that, by transforming a complicated, possibly time-varying and/or path-dependent point process to a unit\edit{-}rate Poisson process, one can more easily observe and quantify exceedances from confidence intervals for the model. Both the simulation and statistical techniques are perhaps most often used for non-stationary Poisson processes, where one can immediately petition to the classical Poissonian idea of conditional uniformity. Here, parking functions will allow us to extend this notion to the clusters within Hawkes processes. 

Our nearest predecessors in random time change theory are likely \citet{giesecke2005dependent} and \citet{ding2009time}. Both of these papers are interested in using this concept to build stochastic models. \citet{giesecke2005dependent} offer\edit{s} somewhat of a converse to the random time change theorem, inverting from unit\edit{-}rate Poisson epochs to construct a more general point process. \citet{ding2009time} uses a similar idea to construct point processes for financial contexts by converting from a pure birth process. By comparison, in this work we will not construct a point process generally, but rather uncover structure specific to the Hawkes process that becomes visible only when the process is transformed akin to the random time change theorem. Furthermore, it is fundamental to our approach that we are inspecting the cluster while conditioning on its size. Our pursuit of the times within a Hawkes cluster is framed by the (conditional) knowledge of how many times there will be, and thus our methodology is indebted conceptually to both the random time change theorem and the conditional uniformity property enjoyed by the Poisson process. This size-conditioning actually constitutes a departure from the typical random time change theorem assumptions in subtle yet important ways, as we will discuss in detail.

This brings us to our contributions and to the organization of the remainder of the paper. Our methodological goal in this work is to use one of the most important results for the Hawkes process, \citet{hawkes1974cluster}'s cluster definition, and one of the most important results for point processes in general, the random time change theorem, to create an analog of one of the most powerful tools for the Poisson process, conditional uniformity (Lemma~\ref{Gentimechange}). This is powered by an unexpected connection between two notable stochastic models, Hawkes process clusters and random parking functions. In linking the continuous time point process to this discrete random vector, we will find that parking functions uncover the full chronology of Hawkes clusters (Theorem~\ref{distThm}). Through this, we are brought back to our original goal, which is to understand the cluster duration. We organize this analysis as follows. In Section~\ref{backgroundSec}, we review the key background concepts, providing precise definitions  and summarizing important results from the literature. Then, in Section~\ref{randTimeSec}, we invoke the compensator transform of the Hawkes process and adapt the random time change theorem for our setting. Our primary result, in which we formalize the Hawkes-process-parking-function connection, is shown in Section~\ref{decompSec}. To demonstrate its usefulness, we apply the techniques to the popular Markovian Hawkes process in Section~\ref{markovSec}, where we use conditional uniformity and parking functions to prove an explicit equivalence between the cluster duration and a random sum of conditionally independent exponential random variables (Theorem~\ref{markovDistThm}). In a second application of the parking function decomposition, in Section~\ref{simSec} we describe the resulting simulation algorithm, which may be of importance to applications in rare event simulation thanks to its ability to explicitly condition on the cluster size (Algorithm~\ref{algSim}). Finally, in Section~\ref{concSec} we conclude.



\section{Model, Scope, and Background Fundamentals}\label{backgroundSec}

To begin, let us briefly review the definitions and prior results that set the stage for our analysis. In particular, we will provide formal definitions of the two focal objects in this work, the Hawkes process cluster and the parking function, while also touching on related topics like branching processes, the Borel distribution, Poisson conditional uniformity, and Dyck paths. 

\subsection{Hawkes Processes, Branching Processes, and Poisson Conditional Uniformity}


Originally introduced by \citet{hawkes1971spectra}, the Hawkes point process $N_t$ and intensity $\lambda_t$ are defined such that
\begin{align}
\PP{N_{t+\delta} - N_t = n \mid \mathcal{F}_t} 
&
=
\begin{cases}
\lambda_t \delta + o(\delta) & n = 1\\
1 - \lambda_t \delta + o(\delta) & n = 0\\
o(\delta) & n > 1
\end{cases}
\label{mainHawkesDef1}
\end{align}
where $\mathcal{F}_t$ is the natural filtration on the underlying probability space $(\Omega, \mathcal{F}, \mathbb{P})$ generated by the point process $N_t$, and where the conditional intensity function $\lambda_t$ is given by
\begin{align}
\lambda_t
&
=
\lambda 
+
\int_{-\infty}^t g(t-u) \mathrm{d}N_u 
=
\lambda
+
\sum_{A_i < t}g(t- A_i)
,
\label{mainHawkesIntensityDef1}
\end{align}
with $\{A_i \mid i \in \mathbb{Z}\}$ as the increasing sequence of arrival epochs in the point process $N_t$. The function $g:\mathbb{R}_+ \to \mathbb{R}_+$ governs the excitement generated upon each arrival epoch, and thus is often referred to as the \textit{excitation kernel} or excitation function. On the other hand, $\lambda \geq 0$ is commonly called the \textit{baseline intensity}, as it drives an underlying stream of exogenous arrivals. By construction in \eqref{mainHawkesIntensityDef1}, every event epoch increases the intensity upon its occurrence by $g(0)$ and then generally $s$ time units later by $g(s)$, hence earning the process its hallmark trait of ``self-excitement.'' For this reason, it is the history of the process that drives its future activity.


One of the most powerful tools for analysis of Hawkes processes has been an alternate definition of the model, originally provided and proved to equivalent to~\eqref{mainHawkesDef1} and~\eqref{mainHawkesIntensityDef1} in \citet{hawkes1974cluster}. This formulation, frequently referred to as the cluster-based definition, bears a style similar to branching process. In a first-level stream, initial events are generated according to a Poisson process at rate $\lambda$. Then, in a secondary stream for each initial event, a progeny cluster is generated independently of all other clusters and arrivals. Starting with the initial event and for each successive arrival in the cluster, direct offspring of that event are generated according to an inhomogeneous Poisson process with rate $g(t)$ for $t$ time units elapsed since the given arrival epoch. This repeats with generating descendants of the offspring themselves until no further arrivals occur in the cluster, and extinction is guaranteed if $\rho = \int_0^\infty g(t) \mathrm{d}t < 1$.

One of the foremost benefits of this alternate definition of the Hawkes process is that it makes very clear the two types of arrivals: the baseline-generated stream and the excitement-generated clusters that spawn off of it. Because the baseline-generated stream is a Poisson process, its behavior is well understood, and so the cluster-based definition allows us to focus on the impact of the self-excitement. This is where the focus of this paper will lie. To dedicate our attention to the structure of the self-excitement, we will narrow  the primary  definition and isolate the clusters. That is, let us essentially mirror Equations~\eqref{mainHawkesDef1} and~\eqref{mainHawkesIntensityDef1} but with a time 0 initial arrival ($A_0 = 0$) and no baseline stream ($\lambda=0$).

\begin{definition}[Hawkes Process Cluster]
\emph{For $t \geq 0$, we will henceforward take $(\lambda_t, N_t)$ as the cluster-specific intensity-and-point-process pair with $N_t$ functioning as in Equation~\eqref{mainHawkesDef1} and with the intensity being simply
\begin{align}
\lambda_t 
= 
\int_0^t g(t-u) \mathrm{d}N_u
=
\sum_{i=0}^{N_t-1} g(t-A_i)
,
\label{mainHawkesDefCluster}
\end{align}
where $A_0 = 0$ without loss of generality. Supposing $\rho = \int_0^\infty g(t) \mathrm{d}t < 1$, the cluster can be fully characterized by its arrival epochs $0 = A_0 < A_1 < \dots < A_{N-1}$, where $N = \lim_{t\to\infty} N_t$ is the cluster size and $\tau = A_{N-1}$ the cluster duration.} 
\end{definition}

One can think of the cluster initializing with a time zero arrival from the exogenous baseline stream, and thus the above definition focuses only on the endogenously driven activity. The stability assumption $\rho < 1$ traces back to Hawkes' original work, and it provides that the size $N$ is finite almost surely and that the intensity $\lim_{t\to\infty} \lambda_t = 0$ almost surely \citep{hawkes1971spectra}. The \citet{hawkes1974cluster} alternate definition of the model already provides \edit{a} perspective of the cluster's descendant structure irrespective of time, and this describes the cluster's size with clarity. Because each event spawns its own offspring Poisson process with time-varying rate $g(t)$, its expected total number of offspring events is $\rho$, and moreover the total number of direct offspring of any one event is $\mathsf{Pois}(\rho)$ distributed. This means that the family tree becomes a Poisson branching process and thus the total progeny will be Borel distributed \citep{feller2008introduction}.  That is, the size $N \sim \mathsf{Borel}(\rho)$ has probability mass function
\begin{align}
\PP{
N
=
k
}
=
\frac{e^{-\rho k}}{k!}\left(\rho k\right)^{k-1}
,
\label{borelDist}
\end{align}
for all $k \in \mathbb{Z}_+$.

While this clear understanding of the cluster size is helpful, on the surface it only tells us a limited part of the story. However, we will see here that \edit{this time-agnostic quantity} provides a key to the chronology as well. Our approach \edit{to unlocking these epochs} is rooted in \edit{the} concept \edit{of} conditional uniformity property for Poisson processes. The classical notion of conditional uniformity states that for a Poisson process with a homogeneous arrival rate, the joint distribution of the epochs $A_1 < A_2 < \dots < A_k$ given that there were $k$ arrivals in the interval $[0,t)$ is equivalent to the joint distribution of the order statistics of $k$ i.i.d.~uniform random variables on $[0,t)$. Alternately stated, as a random vector, $[A_1, A_2, \dots, A_k]$ is conditionally uniform on the polytope $\mathcal{U} = \{x \in [0,t)^k \mid x_i < x_{i+1} ~\forall~ i \leq k - 1\}$. The volume of $\mathcal{U}$ is $t^k \slash k!$, and hence the joint density of this random vector is $k! \slash t^k$ for all $x \in \mathcal{U}$, as the generalized uniform distribution on a polytope means that all points in that polytope have density given by the inverse of the volume. It is of course straightforward to sample from $\mathcal{U}$, as one can simply generate $k$ i.i.d.~$\mathsf{Uni}(0,1)$ random variables, sort them, and multiply by $t$. This creates a handy way of sampling Poisson processes, and the idea extends to time inhomogeneous Poisson processes as well. For a Poisson process with time-varying arrival rate given by $f:\mathbb{R} \to \mathbb{R}^+$ with $F(t) = \int_0^t f(s) \mathrm{d}s$, it can be simulated on an interval $[0,t)$ where $f(\cdot) > 0$ by generating the number of epochs according to $\mathsf{Pois}(F(t))$, sampling the yielded number of standard uniform random variables, sorting them, and returning each arrival epoch as $A_i = F^{-1}(U_{(i)} F(t))$, where $U_{(i)}$ is the $i^\text{th}$ smallest. This idea can be seen to follow simply from the inverse transform sampling method and the likelihood function of the time-varying Poisson process, but it can also be seen as a consequence of the random time change theorem, as we will discuss in Section~\ref{randTimeSec}.


\subsection{Dyck Paths and Parking Functions}

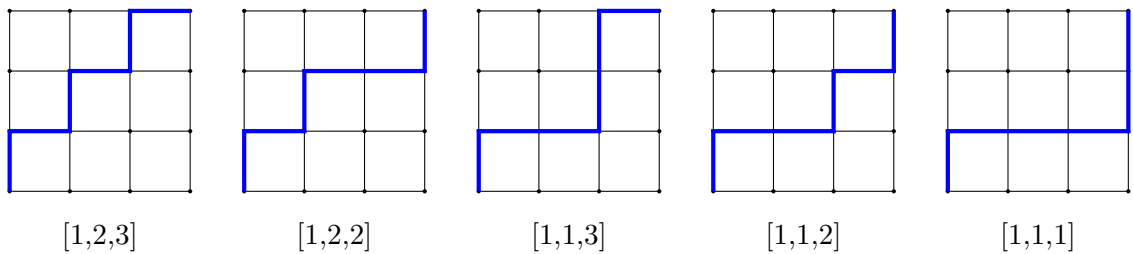
\begin{figure}[h]
\centering
~
~
\hfill
~
\begin{tikzpicture}[scale = 0.8]
\draw (0,0) grid (3,3);
\foreach \i in {0,...,3}
    \foreach \j in {0,...,3}
        \fill (\i,\j) circle (1pt);
\draw[ultra thick,blue] (0,0) |- (1,1) |- (2,2) |- (3,3);
\draw (1.5,-0.75) node {[1,2,3]};
\end{tikzpicture}
~
\hfill
~
\begin{tikzpicture}[scale = 0.8]
\draw (0,0) grid (3,3);
\foreach \i in {0,...,3}
    \foreach \j in {0,...,3}
        \fill (\i,\j) circle (1pt);
\draw[ultra thick,blue] (0,0) |- (1,1) |- (2,2) -| (3,3);
\draw (1.5,-0.75) node {[1,2,2]};
\end{tikzpicture}
~
\hfill
~
\begin{tikzpicture}[scale = 0.8]
\draw (0,0) grid (3,3);
\foreach \i in {0,...,3}
    \foreach \j in {0,...,3}
        \fill (\i,\j) circle (1pt);
\draw[ultra thick,blue] (0,0) |- (1,1) |- (2,1) |- (3,3);
\draw (1.5,-0.75) node {[1,1,3]};
\end{tikzpicture}
~
\hfill
~
\begin{tikzpicture}[scale = 0.8]
\draw (0,0) grid (3,3);
\foreach \i in {0,...,3}
    \foreach \j in {0,...,3}
        \fill (\i,\j) circle (1pt);
\draw[ultra thick,blue] (0,0) |- (1,1) |- (2,1) |- (3,2) -- (3,3);
\draw (1.5,-0.75) node {[1,1,2]};
\end{tikzpicture}
~
\hfill
~
\begin{tikzpicture}[scale = 0.8]
\draw (0,0) grid (3,3);
\foreach \i in {0,...,3}
    \foreach \j in {0,...,3}
        \fill (\i,\j) circle (1pt);
\draw[ultra thick,blue] (0,0) |- (1,1) |- (2,1) |- (3,1) -- (3,3);
\draw (1.5,-0.75) node {[1,1,1]};
\end{tikzpicture}
~
\hfill
~
~
\caption{The collection of all 3-step Dyck paths.}\label{dyckfig}
\end{figure}

Well-known in combinatorics, a $k$-step Dyck path is a non-decreasing trail of points on the lattice from $(0,1)$ to $(k,k)$ such that the height of the path never exceeds the \edit{right-most value of its $x$-coordinate integer interval}. \edit{Here, ``$k$-step'' means that there are $k$ horizontal moves across the lattice in total.} For reference and intuition's sake, Figure~\ref{dyckfig} shows all 3-step Dyck paths. Dyck paths have myriad connections to other combinatorial objects, and many of these connections follow from the fact that they are enumerated by the Catalan numbers \citep[e.g.][]{stanley1999enumerative}. We will let $\mathsf{D}_k$ be the set of all $k$-step Dyck paths. For convenience downstream, our convention will be to record the Dyck path as the vector containing the integer values that are the \edit{largest} $y$-coordinate \edit{attained at each} $x$-coordinate in the lattice path, omitting the terminus. That is, the $2$-step path $(0,1)\to(1,1)\edit{\to(2,1)}\to(2,2)$ would be recorded $[1, 1]$ and likewise $(0,1)\edit{\to(1,1)}\to(1,2)\to(2,2)$ would be $[1,2]$.

\vspace{.1in}
\vspace{-.2in}

Closely related to Dyck paths are parking functions,  which will serve as a focal point throughout this paper. These can be defined through a parsimonious condition.

\begin{definition}[Parking Function]
\emph{For $k \in \mathbb{Z}_+$, $\pi \in \mathbb{Z}_+^k$ is a parking function of length $k$ if and only if it is such that, when sorted such that $\pi_{(1)} \leq \dots \leq \pi_{(k)}$, $\pi_{(i)} \leq i$ for each $i \in \{1, \dots, k\}$.
}
\end{definition}

 Parking functions earn their name from the following intuitive context. Suppose $k$ cars arrive in successive fashion to a strip of $k$ parking spaces, labeled 1 through $k$. Each car $i$ has a preferred space $\pi_i$. If spot $\pi_i$ is available, then the car will park there, otherwise they will take the next available space after $\pi_i$, meaning greater than $\pi_i$. Hence, the constraints $\pi_{(i)} \leq i ~\forall~ 1 \leq i \leq k$ ensure that all cars will be able to park.
For our context, perhaps the most valuable property of parking functions will be that they can be viewed as labeled Dyck paths, or, equivalently, that a sorted parking function is a Dyck path \citep[see, e.g., Section 13.6 of][]{yan2015parking}. That is, every $k$-step Dyck path is also a parking function of length $k$, \edit{and} every parking function of length $k$ can be seen to be a permutation of a $k$-step Dyck path. We will let $\mathsf{PF}_k$ for each $k \in \mathbb{Z}_+$ be the set of all length-$k$ parking functions, and it will be valuable for us to both enumerate these and review an elegant proof of this enumeration. The proof is due to Pollak, but it was recorded and published by contemporaries \citep[e.g.][]{riordan1969ballots,foata1974mappings}, and it is also available on p.~836 in \citet{yan2015parking}.

\begin{lemma}\label{pollak}
There are $(k+1)^{k-1}$ parking functions of length $k \in \mathbb{Z}_+$.
\end{lemma}
\proof{Proof (Pollak's Circle Argument).} Suppose that there are instead $k+1$ parking spaces, and let $\tilde{\pi} \in \{1, \dots, k+1\}^k$ be any $k$-length preference vector for these spaces. Let the cars progress one-by-one so that the $i^\text{th}$ car is the $i^\text{th}$ to pick and suppose, as before, that each car attempts to park in their preferred spot $\tilde{\pi}_i$ and, if unavailable, takes the next space open after it. However, suppose now that the spaces are arranged in a circle, so that if a car has made it to space $k+1$ without parking, it will start again at space 1. Hence, one can instead think of $\tilde{\pi}$ as an element of $(\mathbb{Z}\slash (k+1)\mathbb{Z})^k$. By the pigeonhole principle, all $k$ cars will be able to park and there will be exactly one space remaining, some $\ell \in (\mathbb{Z} \slash (k+1) \mathbb{Z})$. If $\ell = k+1$, then $\tilde \pi$ is a parking function, and if not, it can be converted to one by defining $\pi = (\tilde \pi - \ell)~\mathsf{mod}~(k + 1)$, effectively rotating the empty spot to align with  space $k+1$. Since there are $k+1$ possible rotations that would convert to the same parking function, the cardinality of the set of parking functions is $1\slash(k+1)$ of the cardinality of the set of preference vectors, which is $(k+1)^k$.
\hfill\Halmos\endproof

\begin{figure}[h]
\centering
\begin{tikzpicture}
\foreach \a in {1,2,3,5,6}{
\draw (-\a*360/6 + 150: 1.25cm) node[rotate= -\a*360/6+60]{\a};
}
\draw (-4*360/6 + 150: 1.25cm) node[rotate= -4*360/6+60]{\ul{4}};
\draw[->,thick] +(95:1.6cm) arc(95:-240:1.6cm);
\draw (90: 2cm) node{1 ($\tilde \pi_3$)};
\draw (30: 2.25cm) node{2 ($\tilde \pi_1$)};
\draw (-30: 2.25cm) node{6 ($\tilde \pi_5$)};
\draw (-90: 2cm) node[opacity=0]{4 ($\pi_4$)};
\draw (-150: 2.25cm) node{5 ($\tilde \pi_2$)};
\draw (-210: 2.25cm) node{5 ($\tilde \pi_4$)};
\draw (0: 3.5cm) node{\large$\boldsymbol{\longrightarrow}$};
\end{tikzpicture}
~
\begin{tikzpicture}
\foreach \a in {1,2,3,4,5}{
\draw (-\a*360/6 + 150: 1.25cm) node[rotate= -\a*360/6+60]{\a};
}
\draw (-6*360/6 + 150: 1.25cm) node[rotate= -6*360/6+60]{\ul{6}};
\draw[->,thick] +(95:1.6cm) arc(95:-240:1.6cm);
\draw (90: 2cm) node{1 ($\pi_2$)};
\draw (30: 2.25cm) node{1 ($\pi_4$)};
\draw (-30: 2.25cm) node{3 ($\pi_3$)};
\draw (-90: 2cm) node{4 ($\pi_1$)};
\draw (-150: 2.25cm) node{2 ($\pi_5$)};
\draw (-210: 2.25cm) node[opacity=0]{5 ($\tilde \pi_4$)};
\end{tikzpicture}
\caption{Example conversion of a preference vector drawn from $\boldsymbol{\{1, \dots, k+1\}^k}$ (here, $\boldsymbol{\tilde \pi = [2, 5, 1, 5, 6]}$) to a parking function of length $\boldsymbol{k}$ ($\boldsymbol{\pi = [4, 1, 3, 1, 2]}$) by rotating the empty space to position $\boldsymbol{k+1}$.}\label{pollakfig}
\end{figure}
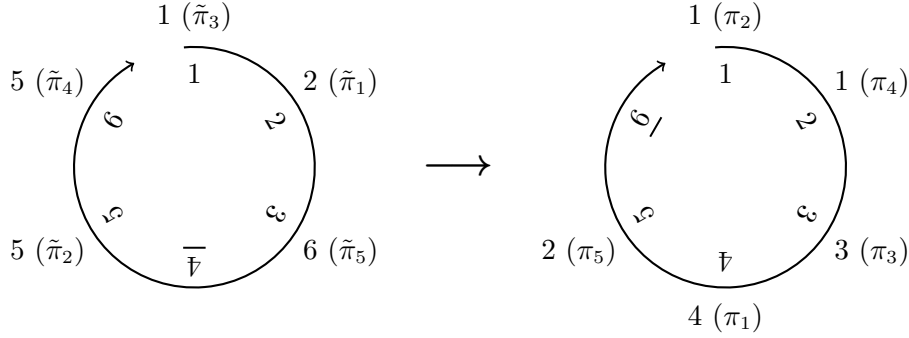

A demonstration of this argument and the concept of rotating the empty space is shown in Figure~\ref{pollakfig}. This construction along the circle leads to a group theoretic sampling of uniformly random parking functions that is used in both~\citet{diaconis2017probabilizing} and~\citet{kenyon2021parking}, and we will employ this idea later for our proposed Hawkes cluster simulation algorithm in Section~\ref{simSec}.

\section{Conditional Uniformity through Random Time Change}\label{randTimeSec}


To connect conditional uniformity and Hawkes processes and to begin understanding this relationship, we first need to define the \textit{compensator} of the process, a transformation of the stochastic process. Let $G(t) = \int_0^t g(u) \mathrm{d}u$. Then, the continuous time compensator of the Hawkes cluster intensity is given by
\begin{align}
\Lambda(t)
=
\int_0^t \lambda_s \mathrm{d}s
=
\sum_{i=0}^{N_t-1}
\int_{A_i}^t
g(s-A_i)
\mathrm{d}s
=
\sum_{i=0}^{N_t-1} G(t-A_i)
.
\end{align}
By construction, $\Lambda(t)$ is a continuous and increasing function. In general, the difference of a point process and its compensator, $N_t - \Lambda(t)$, is a martingale, which can be seen from the Doob-Meyer decomposition \citep[see, e.g., Lemma 7.2.V in][]{daley2003introduction}, and this hints at some of the elegant martingale approaches to the random time change theorem.
In the case of the Hawkes process, we can see that the process is deterministic between its arrival epochs, and so we can capture the full information of the sample path through the discrete time compensator process evaluated only at these points in time:
\begin{align}
\Lambda_k
=
\Lambda(A_k)
=
\sum_{i=0}^{k-1} G(A_k-A_i)
,
\label{disCompDef}
\end{align}
for $k \in \mathbb{Z}_+$, where $\Lambda_k$ is only defined for clusters of size at least $k+1$, so that $A_k < \infty$. Here, we can see that not only is the sequence increasing by definition since $\lambda_t > 0$ \edit{for any finite $t$}, we can also recall that $\lim_{t\to\infty} G(t) = \rho =  \int_0^\infty g(u) \mathrm{d}u$, and so
$
\Lambda_{k-1} < \Lambda_k < k\rho.
$
The compensator sets the stage for us to invoke and extend the classical random time change theorem.

\begin{lemma}\label{Gentimechange}
Given that there are $N = k + 1$ events in the cluster in total including the initial time 0 event, the conditional joint density of the transformed epochs is
\begin{align}
f(\Lambda_1, \Lambda_2, \dots, \Lambda_k \mid N = k + 1)
=
\left(\frac{1}{\rho}\right)^k\frac{(k+1)!}{(k+1)^k}
,
\end{align}
for all $k \in \mathbb{N}$ and all $0 < \Lambda_1 < \dots < \Lambda_k < k\rho$ with $\Lambda_i < i \rho$ for all $i \in \{1, \dots, k\}$. Hence, $\Lambda \mid N \sim \mathsf{Uni}(\mathcal{P}_{N-1})$ for $\Lambda$ as the vector of the compensator points and the polytope ${\mathcal{P}}_k = \{x \in \mathbb{R}^k_+ \mid x_i < i \rho ~\forall~ 1 \leq i \leq k , x_i < x_{i+1} ~\forall~ 1 \leq i \leq k-1\}$.
\end{lemma}
\proof{Proof.}
From the classical random time change theorem, e.g.~Equation (2.14) of \citet{brown2002time}, the unconditioned joint density of the first $k$ compensator points is 
\begin{align*}
f(\Lambda_1, \Lambda_2, \dots, \Lambda_k) 
&=
e^{-\Lambda_k}
,
\end{align*}
where $\Lambda_i < i\rho$ for all $1 \leq i \leq k$ is both implied and required whenever there are at least $k$ offspring events. Then, the conditional density can be expressed
\begin{align*}
f(\Lambda_1, \Lambda_2, \dots, \Lambda_k \mid N = k + 1)
&=
\PP{N = k+1 \mid \Lambda_1, \Lambda_2, \dots, \Lambda_k}
\frac{
e^{-\Lambda_k}
}{
\PP{N = k+1}
}
.
\end{align*}
Now, given that there have been $k$ epochs (excluding the initial time 0 \edit{arrival}) so far up to time $A_k$, the event that there will be exactly $k$ epochs in the cluster in total is equivalent to the event that there will be no epochs after time $A_k$, i.e.~$N_{[A_k,\infty)} = N - N_{A_k} = 0$. Given the history of the process up to time $A_k$ (as contained in the natural filtration $\mathcal{F}_{A_k}$), the probability that $N_{[A_k,\infty)} = 0$ is given by
\begin{align*}
\PP{N_{[A_k,\infty)} = 0 \mid \mathcal{F}_{A_k}}
&=
e^{-\int_{A_k}^\infty \lambda_t \mathrm{d}t}
=
e^{-\sum_{i=0}^{k}\int_{A_k}^\infty g(t -A_i) \mathrm{d}t}
,
\end{align*}
since Equation~\eqref{mainHawkesDef1} shows that the Hawkes process is conditionally non-stationary Poisson given its history. Because the intensity $\lambda_t$ must be strictly positive on $[0,A_k]$, there is a unique and deterministic bijection between $(A_1, A_2, \dots, A_k)$ and $(\Lambda_1, \Lambda_2, \dots, \Lambda_k)$. Hence, because the collection of arrival epochs completely specifies the sample path of the process, so does the collection of compensator points. Therefore, we have
\begin{align*}
\PP{N = k+1 \mid \Lambda_1, \Lambda_2, \dots, \Lambda_k}
&=
\PP{N_{[A_k,\infty)} = 0 \mid \mathcal{F}_{A_k}}
=
e^{-\sum_{i=0}^{k}\int_{A_k}^\infty g(t -A_i) \mathrm{d}t}
,
\end{align*}
which by the definition of the compensator further implies that
\begin{align*}
\PP{N = k+1 \mid \Lambda_1, \Lambda_2, \dots, \Lambda_k} e^{-\Lambda_k}
=
e^{-\sum_{i=0}^{k}\int_{A_k}^\infty g(t -A_i) \mathrm{d}t -\sum_{i=0}^{k}\int_{0}^{A_k} g(t -A_i) \mathrm{d}t}
=
e^{-(k+1)\rho}
.
\end{align*}
Recalling that $N \sim \mathsf{Borel}(\rho)$, we can substitute this probability mass function for $\PP{N=k+1}$ and achieve the stated result. 
\hfill \Halmos \endproof

Let us briefly contrast Lemma~\ref{Gentimechange} with the classical version of the random time change theorem. It is well-established that the seminal result holds for any simple, adapted point process, and that of course includes the Hawkes process we study here. So, it is not that we are newly extending to self-exciting processes; rather, the novelty is in the use of conditioning. Specifically, in this lemma we obtain the conditional joint density of the compensator points given the size of the cluster.  Because this conditioning is on a tally collected at the end of time, we depart from some of the common assumptions of the classical random time change theorem.  That is, the statement of the theorem typically contains an assumption that $\Lambda(t)$ is not bounded almost surely \citep[such as in Theorem 7.4.I in][]{daley2003introduction}, but here we have discussed that each $\Lambda_i$ is strictly less than $i\rho$. In fact, $\PP{\lim_{t\to\infty} \Lambda(t) = \edit{(k+1)} \rho \mid N=k+1} = 1$.\footnote{\edit{This observation can actually be used to construct an alternate proof of Lemma~\ref{Gentimechange} through a limit of the results in \citet{brown2002time}, arising out of Equations (2.5) and (2.14) therein.}} By direct consequence, we also do not require an infinite sequence of epochs on the positive real half-line \citep[as in Proposition 7.4.IV or Theorem T16 in, respectively,][]{daley2003introduction,bremaud1981point}. Similarly, there is also often an assumption that $\lambda_t > 0$ for all $t$ in the time interval of the transform \citep[as in][]{brown2002time}, but here we are explicitly assuming that the intensity converges to 0: $\PP{\lim_{t\to\infty} \lambda_t = 0 \mid N=k+1} = 1$.

These changes are subtle in assumption but important in consequence, as it shifts the connection's end result from the Poisson to the uniform distribution instead. The key takeaway from Lemma~\ref{Gentimechange} is that the conditional density is constant: it has no dependence on the values of the compensator points. Thus, for a cluster with $k$ post-initial events ($N = k +1$), any collection of compensator points satisfying $0 < \Lambda_1 < \Lambda_2 < \dots < \Lambda_k < k\rho$ is equally likely to occur. As used in the proof of Lemma~\ref{Gentimechange}, the compensator vector entirely characterizes the cluster, so we now have access to the full sample path through these conditionally uniform points. This also shows that any two Hawkes processes with the same $\rho$ will have equivalent distributions of  compensator points; neither the size of the cluster nor the conditional joint density of the compensators depend on $g(\cdot)$ outside of $\rho$.

Still, the establishment of conditional uniformity is not out of the blue. Setting aside the departure from typical random time change theorem assumptions, the conditional uniformity in Lemma~\ref{Gentimechange} is not overly surprising given the well known transformation to a unit-rate Poisson process. However, what is special for this result is the polytope on which the transformed points lie uniformly distributed. Recalling the review of Poisson conditional uniformity in Section~\ref{backgroundSec}, we know that the only constraints on the Poisson arrival epochs are that each point is less than the next. However, Lemma~\ref{Gentimechange} shows that in addition to that, each compensator point in the Hawkes cluster is constrained to be less than the product of its index and $\rho$. Hence, the interesting piece of this conditional uniformity that remains to be understood is determining what exactly it means to be uniformly random on this convex polytope $\mathcal{P}_k$.



\section{Decomposing Clusters through Conditional Uniformity}\label{decompSec}


This leads us to our first main result. Lemma~\ref{Gentimechange} gave us the family of polytopes $\{\mathcal{P}_k \mid k \in \mathbb{N}\}$ upon which the cluster's compensator points are conditionally uniform.  What we will find now is that this implies that parking functions provide a hidden spine \edit{in} these transformed arrival epochs. Theorem~\ref{distThm} shows that uniformly random parking functions actually provide a partition for these polytopes from which the compensators can be drawn directly through standard uniforms. Hence, this decomposition elucidates the structure within self-excitement and unites Hawkes processes with parking functions.

\begin{theorem}\label{distThm}
Let $k \in \mathbb{N}$. Given that the Hawkes process cluster is comprised of $N=k+1$ events, the compensator transformed arrival epochs $\Lambda_1, \dots, \Lambda_k$ are such that
\begin{align}
\Lambda_i 
\stackrel{\mathsf{D}}{=}
\rho\cdot\left(\pi - U\right)_{(i)}
\qquad
\forall\,i \in \{1, \dots, k\}
,
\end{align}
where $\left(\pi - U\right)_{(1)} < \dots < \left(\pi - U\right)_{(k)}$, with $\pi \in \mathsf{PF}_k$ as a $k$-step parking function drawn uniformly at random and $U \in \edit{(}0,1)^k$ as $k$ i.i.d.~standard uniform random variables.
\end{theorem}




We will prove Theorem~\ref{distThm} across a series of subsidiary statements. For simplicity and without loss of generality, let us consider the unit polytope $\bar{\mathcal{P}}_k = \{x \in \mathbb{R}^k_+ \mid  x_i < i ~\forall~ 1 \leq i \leq k, x_i < x_{i+1} ~\forall~ 1 \leq i \leq k-1\}$, which maps to $\mathcal{P}_k$ through the bijection $x \mapsto \rho x$ for all $k$. To begin to establish our methodology of sampling from $\bar{\mathcal{P}}_k$, we will first view $\bar{\mathcal{P}}_k$ as a union over sub-polytopes with measure 0 intersections. Specifically, let us distill each point in $\bar{\mathcal{P}}_k$ as lying on intervals bounded by integer steps in each coordinate direction. Then, any point in $\bar{\mathcal{P}}_k$ can be seen to lie in one of these subspaces. 

To organize this decomposition of $\bar{\mathcal{P}}_k$, we will draw in the first combinatorial object discussed in Section~\ref{backgroundSec}: Dyck paths. In Proposition~\ref{pDecompProp}, we use the collection of $k$-step Dyck paths to partition $\bar{\mathcal{P}}_k$ into $|\mathsf{D}_k| = C_k = {2k \choose k}\slash(k+1)$ disjoint subspaces. As the brief proof demonstrates, the result itself is straightforward, but it will be of use in structuring the uniform sampling from the polytope.

\begin{proposition}\label{pDecompProp}
For each $k \in \mathbb{Z}_+$, the unit polytope $\bar{\mathcal{P}}_k$ can be decomposed through the set of $k$-step Dyck paths, i.e.
\begin{align}
\bar{\mathcal{P}}_k
&=
\bigsqcup_{d \in \mathsf{D}_k}
\bar{\mathcal{P}}_{k,d}
,
\label{pDecompEq}
\end{align}
where
\begin{align}
\bar{\mathcal{P}}_{k,d}
&=
\{x \in \mathbb{R}_+^k 
\mid 
d_i-1 \leq x_i < d_i ~\forall~ 1 \leq i \leq k, x_i < x_{i+1} ~\forall~ 1 \leq i \leq k-1
\}
\label{PkdDefEq}
,
\end{align}
for each Dyck path $d \in \mathsf{D}_k$.
\end{proposition}
\proof{Proof.}
By definition, the sets $\bar{\mathcal{P}}_{k,d}$ are mutually disjoint across all paths $d$. Because the ordering constraint $\edit{x_i < x_{i+1}}$ is shared by every set in the Dyck path union and by $\bar{\mathcal{P}}_k$, we only must check that the space constraints agree. Here, we can quickly observe that $1 \leq d_i \leq i$ for each $i \in \{1, \dots, k\}$, hence $\bigcup_{d \in \mathcal{D}_k}
\{x_i \in \mathbb{R}_+
\mid 
d_i-1 \leq x_i < d_i \} =
\{x_i \in \mathbb{R}_+
\mid 
\edit{0 \leq x_i} < i \}$.
\hfill \Halmos \endproof


\begin{figure}
\centering
\includegraphics[width=0.95\textwidth]{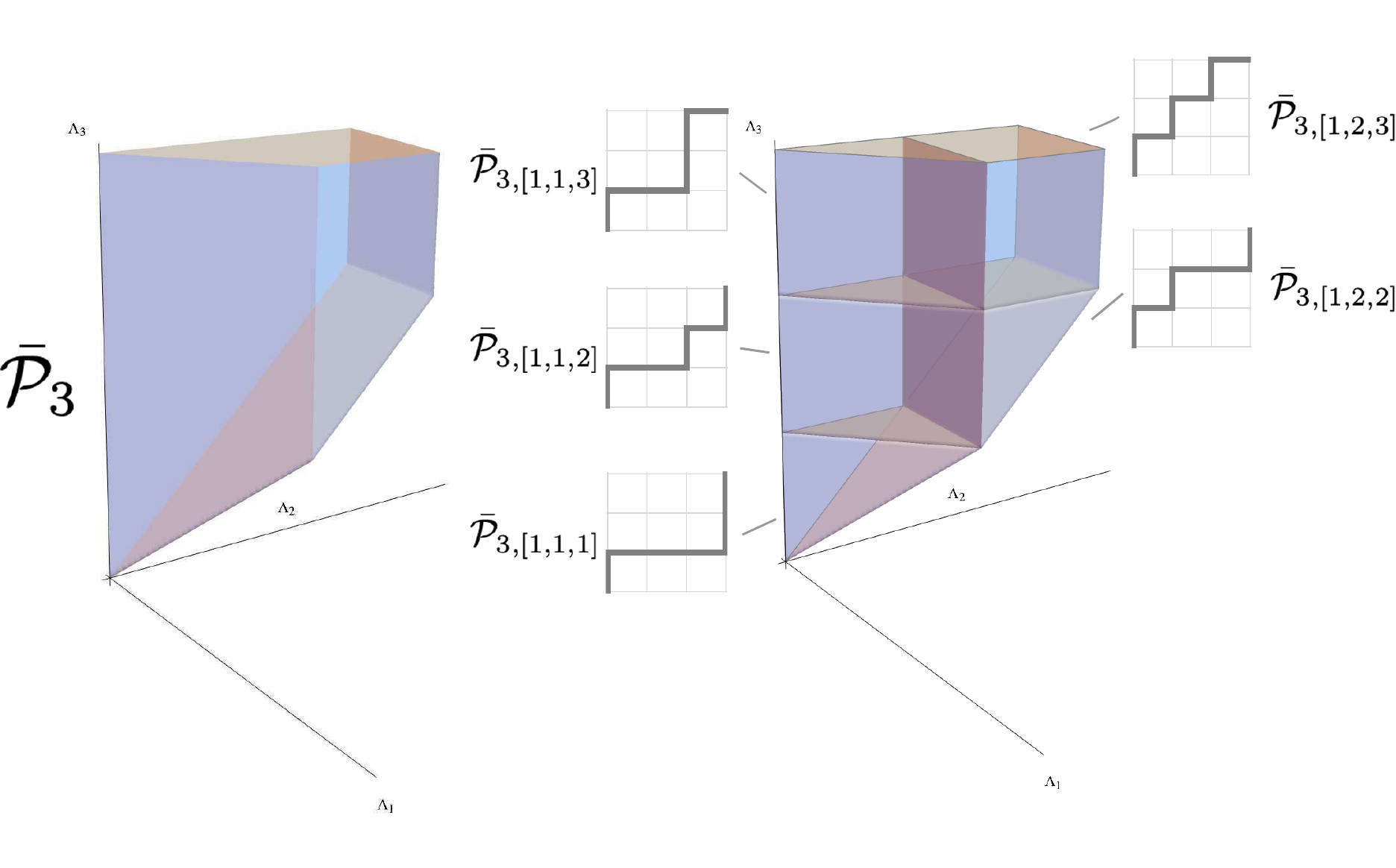}
\caption{The three dimensional polytope (left) and its decomposition into the five sub-regions indexed by Dyck paths (right).}\label{decompFig}
\end{figure}

The three dimensional polytope, $\bar{\mathcal{P}}_3$, and its Dyck path decomposition are shown in Figure~\ref{decompFig}. In some sets in the Dyck path partition, the ordering constraints will be superfluous given the integer ranges on which each coordinate lives. For example, at any dimension $k$ the largest volume subset will be the hypercube $\bar{\mathcal{P}}_{k,[1, 2, \dots, k]} = \{x \in \mathbb{R}_+^k 
\mid 
i-1 \leq x_i < i ~\forall~ 1 \leq i \leq k
\}$, in which the ordering constraints are not needed. By comparison, the smallest volume subset will correspond to the ``lowest'' Dyck path $[1, 1, \dots, 1]$, which yields $\bar{\mathcal{P}}_{k,[1, 1, \dots, 1]} = \{x \in \edit{[0,1)}^k 
\mid 
x_i < x_{i+1} ~\forall~ 1 \leq i \leq k-1
\}$.\footnote{\edit{In fact, $\bar{\mathcal{P}}_{k,[1, 1, \dots, 1]}$ is equivalent (up to scaling) to the polytope on which homogeneous Poisson arrival epochs lie, $\mathcal{U}$. The self-exciting patterns of the Hawkes create the added regions and their corresponding constraints.}} This structure is also intimately linked to the value of the decomposition. The benefit of Proposition~\ref{pDecompProp} is that it is simple and intuitive to sample uniformly within a given subspace from the Dyck path partition. That is, a uniform sample from the hypercube is simply  $U_i \sim \mathsf{Uni}(i-1,i)$, drawn independently for each $i$. On the other hand, a uniform sample from the $[1,1,\dots, 1]$ subspace is given by the order statistics of $k$ independent $\mathsf{Uni}(0,1)$ random variables. This generalizes cleanly. To sample uniformly from the region given by the Dyck path $d \in \mathsf{D}_k$, draw $k$ independent $\mathsf{Uni}(0,1)$ random variables, say $U_1, \dots, U_k$, and return the order statistics of $d_1 - U_1, \dots, d_k - U_k$ as the sample. This is again straightforward to show.

\begin{proposition}\label{subsampleprop}
For $k \in \mathbb{Z}_+$ and $d \in \mathsf{D}_k$, let $X_d \in \mathbb{R}_+^k$ be equivalent in distribution to the order statistics of $d - U$, where $U \in \edit{(}0,1)^k$ is such that $U_i \stackrel{\mathsf{iid}}{\sim} \mathsf{Uni}(0,1)$. Then, $X_d \sim \mathsf{Uni}(\bar{\mathcal{P}}_{k,d})$.
\end{proposition}
\proof{Proof.}
The two distributions share a sample space by definition. Because the elements of $U$ are independent, distinct values of $d$ yield independent elements of $X_d$. When values are repeated in $d$, the corresponding elements of $X_d$ are equivalent to shifted standard uniform order statistics. Hence, for $\kappa_i(d) = |\{j \mid d_j = i\}|$, the density of $X_d$ is $f(x_1, \dots, x_k) = \prod_{i=1}^k \kappa_i(d)!$, which is precisely the inverse of the volume of $\bar{\mathcal{P}}_{k,d}$.
\hfill \Halmos \endproof

Hence, we see that the Dyck path indexing provides a way of encoding which points share an interval and which do not. \edit{To use} this decomposition to uniformly sample from $\bar{\mathcal{P}}_k$, \edit{it remains} to find a distribution over the Dyck paths, so that we draw from $\bar{\mathcal{P}}_k$ by first selecting a set from the partition and then draw a point from within it. Because Dyck paths are enumerated by the Catalan numbers, we know that there will be $C_k = {2k \choose k}\slash(k+1)$ sets within the partition of $\bar{\mathcal{P}}_k$. However, each set should not be equally likely; rather, the likelihood of each set should be proportional to its volume. This is clear at $k=2$. $\bar{\mathcal{P}}_{2,[1, 2]}$, a square comprised of $0 \leq x_1 < 1$ and $1 \leq x_2 < 2$, should be twice as likely as $\bar{\mathcal{P}}_{2,[1, 1]}$, the half-square given by $0 \leq x_1 <
x_2 < 1$. At $k=3$, there is one cube ($\bar{\mathcal{P}}_{3,[1, 2,3]}$), three half-cubes ($\bar{\mathcal{P}}_{3,[1, 1,2]}$, $\bar{\mathcal{P}}_{3,[1, 1,3]}$, and $\bar{\mathcal{P}}_{3,[1, 2,2]}$), and one sixth-cube ($\bar{\mathcal{P}}_{3,[1, 1,1]}$), as shown in Figure~\ref{decompFig}. Since the total volume of $\bar{\mathcal{P}}_3$ is $8/3$, the respective probabilities for each of the five regions should be $3/8$, $3/16$, $3/16$, $3/16$, and $1/16$.


In what is of beautiful consequence for our problem, parking functions yield precisely the correct weighting of Dyck paths. Recalling from Section~\ref{backgroundSec} that there is a bijection between sorted parking functions and Dyck paths, the suite of proportions of length-$k$ parking functions that \edit{sort} to each $k$-step Dyck path  yields exactly the desired distribution over the partitions of $\bar{\mathcal{P}}_k$. Hence, to properly select the partition subset from which to sample in $\bar{\mathcal{P}}_k$, we can simply sample the Dyck paths according to the fractions of equivalent parking functions, or, even more directly, we can merely draw a parking function uniformly at random, sort it, and select the corresponding subset of the polytope. We prove this now in Proposition~\ref{PFdistProp}.

\begin{proposition}\label{PFdistProp}
For $k \in \mathbb{Z}_+$, let $\pi \in \mathsf{PF}_k$ be a uniformly random parking function of length $k$, and let $\vec{\pi}$ be the parking function sorted increasingly. Likewise, let $\Lambda \in \mathbb{R}_+^k$ be selected uniformly at random from $\bar{\mathcal{P}}_k$. Then, for any Dyck path $d \in \mathsf{D}_k$,
\begin{align}
\PP{\Lambda \in \bar{\mathcal{P}}_{k,d}}
=
\PP{\vec\pi = d}
=
\frac{k!}{(k+1)^{k-1}
\prod_{i=1}^k \kappa_i(d)!}
,
\label{PFdistEq}
\end{align}
where $\kappa_i(d) = |\{j \mid d_j = i\}|$ counts the number of occurrences of $i$ within a Dyck path $d$.
\end{proposition}
\proof{Proof.}
Let $d \in \mathsf{D}_k$ be arbitrary. We will prove the result by showing that $\PP{\Lambda \in \bar{\mathcal{P}}_{k,d}}$ and $\PP{\vec\pi = d}$ are both equal to the right-most side of Equation~\eqref{PFdistEq}. Starting with $\PP{\vec\pi = d}$, let us consider the space of parking functions with length $k$. From Lemma~\ref{pollak}, there are $(k+1)^{k-1}$ parking functions of length $k$. Then, accounting for repeated values, there are $k!\slash (\prod_{i=1}^k \kappa_i(d)!)$ ways to permute the $k$ values in $d$, with each one being a possible draw from the space of parking functions. Hence, Equation~\eqref{PFdistEq} shows the true distribution of $\vec \pi$. 

Now, for $\PP{\Lambda \in \bar{\mathcal{P}}_{k,d}}$, let us recall that Lemma~\ref{Gentimechange} implies that the density of $\Lambda$ on $\bar{\mathcal{P}}_k$ is $k! \slash (k+1)^{k-1}$. Hence, the probability that $\Lambda$ lies in $d$'s partition set is
\begin{align}
\PP{\Lambda \in \bar{\mathcal{P}}_{k,d}}
&=
\int_0^1\int_{(d_2-1 \vee x_1)}^{d_2} \dots \int_{(d_k-1 \vee x_{k-1})}^{d_k} 
\frac{k!}{(k+1)^{k-1}} 
\mathrm{d}x_k  
\dots
\mathrm{d}x_2
\mathrm{d}x_1
.
\label{intEqP3}
\end{align}
If $d_{i-1} < d_i < d_{i+1}$, \edit{then neither} the range of integration nor the integrand of the integral with respect to $x_i$ will depend on any of the other variables of integration, hence this integration is separable. Moreover, the integrals over $x_i$ and $x_j$ are not separable if and only if $d_i = d_j$. \edit{That is, for any Dyck path $d$, this separability allows us to decompose the integral in Equation~\eqref{intEqP3} into a product over smaller integrals, each of which contains only the variables that share a corresponding $d$ value.} \edit{Ignoring the $k!\slash(k+1)^{k-1}$ coefficient temporarily for brevity's sake}, we can re-express the integral  as
\begin{align*}
\int_0^1\int_{(d_2-1 \vee x_1)}^{d_2} \dots \int_{(d_k-1 \vee x_{k-1})}^{d_k} 
1 \cdot
\mathrm{d}x_k  
\dots
\mathrm{d}x_2
\mathrm{d}x_1
&=
\prod_{i=1}^k
\int_{\edit{i-1}}^{\edit{i}}
\int_{\xi_1}^{\edit{i}}
\dots
\int_{\xi_{\kappa_i(d)-1}}^{\edit{i}}
1 \cdot
\mathrm{d}\xi_{\kappa_i(d)}
\dots
\mathrm{d}\xi_2
\mathrm{d}\xi_1
,
\end{align*}
where an empty integral is taken to equal 1 by convention. In this simple form we can quickly see that
\begin{align*}
\prod_{i=1}^k
\int_{\edit{i-1}}^{\edit{i}}
\int_{\xi_1}^{\edit{i}}
\dots
\int_{\xi_{\kappa_i(d)-1}}^{\edit{i}}
1 \cdot
\mathrm{d}\xi_{\kappa_i(d)}
\dots
\mathrm{d}\xi_2
\mathrm{d}\xi_1
=
\prod_{i=1}^k \frac{1}{\kappa_i(d)!}
,
\end{align*}
which yields the stated result for $\PP{\Lambda \in \bar{\mathcal{P}}_{k,d}}$.
\hfill \Halmos \endproof


Theorem~\ref{distThm} now follows directly.

\proof{Proof of Theorem~\ref{distThm}.}
By Lemma~\ref{Gentimechange}, we have that \edit{the} vector of transformed epochs $\Lambda = [\Lambda_1, \dots, \Lambda_k]$ is jointly uniform over the polytope \edit{$\mathcal{P}_k$}, i.e.~$\Lambda \sim \mathsf{Uni}(\mathcal{P}_k)$. Equivalently, $\Lambda \slash \rho \sim \mathsf{Uni}(\bar{\mathcal{P}}_k)$. Then, by the decomposition into disjoint subsets in Proposition~\ref{pDecompProp}, for $d \in \mathsf{D}_k$ we further have that $\Lambda \slash \rho \mid (\Lambda \slash \rho \in \bar{\mathcal{P}}_{k,d}) \sim \mathsf{Uni}(\bar{\mathcal{P}}_{k,d})$. By Proposition~\ref{subsampleprop}, the order statistics of $d-U$ are also uniform within $\bar{\mathcal{P}}_{k,d}$. Now, Proposition~\ref{PFdistProp} establishes that the probability $\PP{\Lambda \slash \rho \in \bar{\mathcal{P}}_{k,d}} = \PP{\vec{\pi} = d}$ for a uniformly random parking function $\pi$, and thus the order statistics of $\pi - U$ and $\Lambda \slash \rho$ are equivalent in their uniform distribution on $\bar{\mathcal{P}}_k$.
\hfill \Halmos \endproof

\section{Application to the Markovian Hawkes Process}\label{markovSec}

In the remainder of this paper, we will demonstrate some of the advantages of the parking function decomposition of the Hawkes process clusters. We start by showcasing the analytical value of the conditional uniformity through application to what is perhaps the most widely used Hawkes excitation kernel. Within this section, let us now assume that $g(x) = \alpha e^{-\beta x}$ for $0 < \alpha < \beta$. Under this kernel, the intensity becomes a Markov process \citep{oakes1975markovian}. Much of the popularity of the exponential kernel follows from its frequent tractability, as we will see here. 

Let us start by expressing the compensator for the exponential kernel. Given $g(x) = \alpha e^{-\beta x}$, we see that $G(x) = \rho(1-e^{-\beta x})$ for $\rho = \alpha \slash \beta$, and thus
\begin{align*}
\Lambda(t)
=
\sum_{i=0}^{N_t-1}
\rho\left(1-e^{-\beta(t-A_i\edit{)}}\right)
\quad
\text{ and }
\quad
\Lambda_k
=
\sum_{i=0}^{k-1}
\rho\left(1-e^{-\beta(A_k - A_i)}\right)
.
\end{align*}
There is a great deal of information hidden in this expression. For example, if we let $\lambda_k = \lambda_{A_k^-} = \sum_{i=0}^{k-1} \alpha e^{-\beta(A_k - A_i)}$ be the value of the intensity immediately before the cluster's $k^\text{th}$ event, then we can observe that this quantity is merely an affine transformation of the corresponding compensator point $\Lambda_k$. That is, 
\begin{align}
\Lambda_k 
=
\rho k
-
\frac{1}{\beta}
\sum_{i=0}^{k-1}
\alpha e^{-\beta(A_k - A_i)}
=
\rho k
-
\frac{1}{\beta}
\lambda_k
\longrightarrow
\lambda_k
=
\alpha k - \beta \Lambda_k 
\label{intensEq}
.
\end{align}
An interesting consequence of this relationship is that, after shifting and normalizing, the pre-event intensity values will satisfy the same distributional properties seen in Lemma~\ref{Gentimechange} and Theorem~\ref{distThm}. Hence, the conditional uniformity property is even stronger for this kernel. Given the number of events in the cluster, here we find the stronger and even more surprising implication that the intensity sample paths are conditionally uniform. \edit{This is an immediate corollary of Lemma~\ref{Gentimechange}.} The connection to parking functions also carries over naturally, as Equation~\eqref{intensEq} implies $\lambda_k \stackrel{\mathsf{D}}{=}\alpha \left(k - (\pi - U)_{(k)}\right)$.

This tractability also allows us direct access to the distribution of the cluster's duration, as the arrival \edit{epochs} and inter-arrival times can be easily obtained from the compensator points. Letting $S_k = A_k - A_{k-1}$ be the inter-arrival times of the cluster, we can see that
\begin{align*}
\Lambda_k 
= 
\sum_{i=0}^{k-1}
\rho\left(1-e^{-\beta(S_k + A_{k-1} - A_i)}\right)
=
k\rho - e^{-\beta S_k} \left( k\rho - \Lambda_{k-1}\right)
,
\end{align*}
and thus each inter-arrival time is given by
\begin{align}
S_k
=
-\frac{1}{\beta}\log\left(
\frac{k\rho - \Lambda_k}{k\rho - \Lambda_{k-1}}
\right)
.
\end{align}
Hence, we also immediately have that the $k^\text{th}$ \edit{arrival epoch} can be written
\begin{align}
A_k
=
-\frac{1}{\beta}\log\left(
\prod_{i=1}^k
\frac{i\rho - \Lambda_i}{i\rho - \Lambda_{i-1}}
\right)
,
\label{EqMarkovA}
\end{align}
for all $k \in \mathbb{Z}^+$. This convenient form enables us to analyze this paper's most sought after single quantity: the cluster duration $\tau = A_{N-1}$. This brings us to our second main result, in which we provide a parsimonious characterization of the duration \edit{distribution}, showing that $\tau$ is equivalent to a sum of conditionally independent exponential random variables. Again, a uniformly random parking function appears as a hidden spine of the Hawkes process cluster.

\begin{theorem}\label{markovDistThm}
Let $N \sim \mathsf{Borel}(\rho)$ for $\rho = \alpha \slash \beta$, and, given $N$, let $\pi$ be a uniformly random parking function of length $N-1$. Then, the duration of a cluster in the Markovian Hawkes process satisfies
\begin{align}
\tau
\stackrel{\mathsf{D}}{=}
\frac{1}{\beta}
\sum_{i=1}^{N-1}
T_{\pi, i}
,
\end{align}
 where $T_{\pi, i} \sim \mathsf{Exp}\big(i+1 - \sum_{j=N-i}^{N-1} \kappa_j(\pi) \big)$ are conditionally independent given $\pi$, with $\kappa_i(\pi) = |\{j \mid \pi_j = i\}|$.
\end{theorem}
\proof{Proof.}
Let us first find an expression for the conditional Laplace-Stieltjes transform (LST) directly from Lemma~\ref{Gentimechange}, and then interpret it through the lens of Theorem~\ref{distThm}. If $N = 1$, then the cluster contains only the initial time-0 event. Thus, $\PP{\tau = 0 \mid N = 1} = 1$ and $\E{e^{-\theta \tau} \mid N = 1} = 1$. Then, for $k \geq 1$, we can see through applying Equation~\eqref{EqMarkovA} that the conditional LST can be found through the integral
\begin{align*}
\E{e^{-\theta\tau} \mid N = k + 1}
=
\edit{\frac{1}{\rho^k}} \frac{(k+1)!}{(k+1)^k}
\int_{0}^{\rho}
\int_{x_1}^{2\rho}
\dots
\int_{x_{k-1}}^{k\rho}
\left(
\prod_{i=1}^k
\frac{i\rho - x_i}{i\rho - x_{i-1}}
\right)^{\frac{\theta}{\beta}}
\mathrm{d}x_k
\dots
\mathrm{d}x_2
\mathrm{d}x_1
.
\end{align*}
Now, for reference, let us note that for $a  > 0$, \edit{$b \geq 0$,} $c \geq 0$,  and $m \in \mathbb{N}$, we can obtain the following \edit{expression for a soon-to-be relevant} definite integral directly from the binomial theorem:
\begin{align}
\int_0^{1+c} x^m \left(\frac{ax}{1+c}\right)^b \mathrm{d}x
&=
\frac{(1+c)^{m+1} a^b }{m+b+1}
=
\sum_{\ell=0}^{m+1} {m+1 \choose \ell} \frac{a^b c^{\ell}}{m+b+1}
\label{IntIndEq}
.
\end{align}
To make use of this for the conditional LST, let us substitute $y_i = i -  x_i\slash{\rho}$ for each index of integration $x_i$. After iterative application of \eqref{IntIndEq}, the exponent on the next index of integration will depend on the exponent of the previous. This leaves us with
\begin{align*}
&
\E{e^{-\theta\tau} \mid N = k + 1}
=
\frac{(k+1)!}{(k+1)^k}
\int_{0}^{1}
\int_{0}^{1+y_1}
\dots
\int_{0}^{1+y_{k-1}}
\left(
\prod_{i=1}^k
\frac{y_i}{1+y_{i-1}}
\right)^{\frac{\theta}{\beta}}
\mathrm{d}y_k
\dots
\mathrm{d}y_2
\mathrm{d}y_1
\\
&=
\frac{(k+1)!}{(k+1)^k}
\int_{0}^{1}
\int_{0}^{1+y_1}
\dots
\int_{0}^{1+y_{k-2}}
\sum_{\ell_1=0}^1
{1 \choose \ell_1}
\frac{y_{k-1}^{\ell_1}}{\frac{\theta}{\beta}+1}
\left(
\prod_{i=1}^{k-1}
\frac{y_i}{1+y_{i-1}}
\right)^{\frac{\theta}{\beta}}
\mathrm{d}y_{k-1}
\dots
\mathrm{d}y_2
\mathrm{d}y_1
\\
&=
\frac{(k+1)!}{(k+1)^k}
\int_{0}^{1}
\int_{0}^{1+y_1}
\dots
\int_{0}^{1+y_{k-3}}
\sum_{\ell_1=0}^1
\frac{{1 \choose \ell_1}}{\frac{\theta}{\beta}+1}
\sum_{\ell_2=0}^{\ell_1+1}
\frac{{\ell_1+1 \choose \ell_2}}{\frac{\theta}{\beta}+\ell_1+1}
y_{k-2}^{\ell_2}
\left(
\prod_{i=1}^{k-2}
\frac{y_i}{1+y_{i-1}}
\right)^{\frac{\theta}{\beta}}
\mathrm{d}y_{k-2}
\dots
\mathrm{d}y_2
\mathrm{d}y_1
\\
&
\,\,\,\vdots
\\
&=
\frac{(k+1)!}{(k+1)^k}
\sum_{\ell_1=0}^1
\frac{{1 \choose \ell_1}}{\frac{\theta}{\beta}+1}
\sum_{\ell_2=0}^{\ell_1+1}
\frac{{\ell_1+1 \choose \ell_2}}{\frac{\theta}{\beta}+\ell_1+1}
\sum_{\ell_3=0}^{\ell_2+1}
\frac{{\ell_2+1 \choose \ell_3}}{\frac{\theta}{\beta}+\ell_2+1}
\dots
\sum_{\ell_{k-1}=0}^{\ell_{k-2}+1}
\frac{{\ell_{k-2}+1 \choose \ell_{k-1}}}{\frac{\theta}{\beta}+\ell_{k-2}+1}
\int_{0}^{1}
y_{1}^{\frac{\theta}{\beta}+\ell_{k-1}}
\mathrm{d}y_1
\\
&=
\frac{(k+1)!}{(k+1)^k}
\sum_{\ell_1=0}^1
\frac{{1 \choose \ell_1}}{\frac{\theta}{\beta}+1}
\sum_{\ell_2=0}^{\ell_1+1}
\frac{{\ell_1+1 \choose \ell_2}}{\frac{\theta}{\beta}+\ell_1+1}
\sum_{\ell_3=0}^{\ell_2+1}
\frac{{\ell_2+1 \choose \ell_3}}{\frac{\theta}{\beta}+\ell_2+1}
\dots
\sum_{\ell_{k-1}\edit{=}0}^{\ell_{k-2}+1}
\frac{{\ell_{k-2}+1 \choose \ell_{k-1}}}{\frac{\theta}{\beta}+\ell_{k-2}+1}
\frac{1}{\frac{\theta}{\beta}+\ell_{k-1}+1}
,
\end{align*}
where now the bounds of summation, rather than the bounds of integration, depend on the magnitude of the prior term. Combining like terms into a product, this can be re-expressed
\begin{align*}
\E{e^{-\theta\tau} \mid N = k + 1}
&=
\frac{(k+1)!}{(k+1)^k}
\sum_{\ell_1=0}^1
{1 \choose \ell_1}
\sum_{\ell_2=0}^{\ell_1+1}
{\ell_1+1 \choose \ell_2}
\dots
\sum_{\ell_{k-1}=0}^{\ell_{k-2}+1}
{\ell_{k-2}+1 \choose \ell_{k-1}}
\frac{1}{\frac{\theta}{\beta}+1}
\prod_{i=1}^{k-1}\frac{1}{\frac{\theta}{\beta}+\ell_{i}+1}
.
\end{align*}
By decomposing the binomial coefficients into the underlying factorial terms, cancelling numerators with preceding denominators, and moving the remaining pieces to the product, we can pull $(k-1)!$ from the leading $(k+1)!$ and recognize a resulting multinomial coefficient. That is, after cancellation, the remaining denominators contain factorials of terms that will sum to $k-1$, and thus with the leading $(k-1)!$ the multinomial arises:
\begin{align*}
&
\E{e^{-\theta\tau} \mid N = k + 1}
=
\frac{(k+1)!}{(k+1)^k}
\sum_{\ell_1=0}^1
\frac{1!}{\ell_1!(1-\ell_1)!}
\dots
\sum_{\ell_{k-1}=0}^{\ell_{k-2}+1}
\frac{(\ell_{k-2}+1)!}{\ell_{k-1}!(\ell_{k-2}+1-\ell_{k-1})!}
\frac{1}{\frac{\theta}{\beta}+1}
\prod_{i=1}^{k-1}\frac{1}{\frac{\theta}{\beta}+\ell_{i}+1}
\\
&=
\frac{(k+1)!}{(k+1)^k}
\sum_{\ell_1=0}^1
\frac{1}{(1-\ell_1)!}
\sum_{\ell_2=0}^{\ell_1+1}
\frac{1}{(\ell_1+1-\ell_2)!}
\dots
\sum_{\ell_{k-1}=0}^{\ell_{k-2}+1}
\frac{1}{\ell_{k-1}!(\ell_{k-2}+1-\ell_{k-1})!}
\frac{1}{\frac{\theta}{\beta}+1}
\prod_{i=1}^{k-1}\frac{\ell_{i-1}+1}{\frac{\theta}{\beta}+\ell_{i}+1}
\\
&=
\frac{k}{(k+1)^{k-1}}
\sum_{\boldsymbol{\ell}\in\mathcal{L}_{k-1}}
{
k-1
\choose
1-\ell_1, \ell_1+1-\ell_2, \dots, \ell_{k-2}+1-\ell_{k-1}, \ell_{k-1}
}
\frac{1}{\frac{\theta}{\beta}+1}
\prod_{i=1}^{k-1}\frac{\ell_{i-1}+1}{\frac{\theta}{\beta}+\ell_{i}+1}
,
\end{align*}
where $\mathcal{L}_{k-1} = \{\boldsymbol{\ell} \in \mathbb{N}^{k-1} \mid 0 \leq \ell_1 \leq 1, 0 \leq \ell_i \leq \ell_{i-1} + 1 \,\,\forall\,\, 2 \leq i \leq k-1\}$.
To simplify further, we can now multiply and divide by $\ell_{k-1}+1$, allowing us to re-index the product and absorb the leading $k$ into the multinomial, yielding a form which we can now interpret:
\begin{align*}
\E{e^{-\theta\tau} \mid N = k + 1}
&=
\frac{1}{(k+1)^{k-1}}
\sum_{\boldsymbol{\ell}\in\mathcal{L}_{k-1}}
{
k
\choose
1-\ell_1, \dots, \ell_{k-2}+1-\ell_{k-1}, \ell_{k-1} + 1
}
\frac{1}{\frac{\theta}{\beta}+1}
\prod_{i=1}^{k-1}\frac{\ell_{i}+1}{\frac{\theta}{\beta}+\ell_{i}+1}
.
\end{align*}

First, let us observe that $\mathcal{L}_{k-1}$ can be viewed as the set of differences between \edit{each given} $k$-step Dyck path and the diagonal $y = x$ line.  \edit{Hence, an isomorphism with $\mathsf{D}_k$} is implied directly from the constraints \edit{that define $\mathcal{L}_{k-1}$}. Following this, the multinomial coefficient can then be viewed as counting the number of length-$k$ parking functions that can be formed using \edit{each given $\ell$, as we can recall that the parking functions of length $k$ are all the possible orderings of the $k$-step Dyck paths \citep[e.g.,][]{yan2015parking}.} That is, the multinomial coefficient counts the number of parking functions with $1-\ell_1$ $k$'s, $\ell_1 +1 - \ell_2$ $(k-1)$'s, and so on. Using the $\kappa_i(\pi)$ notation, we can then re-express the conditional Laplace-Stieltjes transform as a sum over all parking functions:
\begin{align*}
\E{e^{-\theta \tau} \mid N = k + 1}
&
=
\frac{1}{(k+1)^{k-1}} 
\sum_{\pi \in \mathsf{PF}_k}
\frac{1}{\frac{\theta}{\beta}+1}
\prod_{i=1}^{k-1}
\frac{
i+1 - \sum_{j=k+1-i}^k \kappa_j(\pi)
}
{
\frac{\theta}{\beta}+ i+1 - \sum_{j=k+1-i}^k \kappa_j(\pi)
}
.
\end{align*}
Here, we have used the pattern above to substitute $i - \sum_{j=k+1-i}^k \kappa_j(\pi)$ for $\ell_i$. Because there are $(k+1)^{k-1}$ parking functions of length $k$, this can also be viewed as an expectation over a uniformly sampled parking function. That is,
\begin{align*}
\frac{1}{(k+1)^{k-1}} 
\sum_{\pi \in \mathsf{PF}_k}
\frac{1}{\frac{\theta}{\beta}+1}
\prod_{i=1}^{k-1}
\frac{
i+1 - \sum_{j=k+1-i}^k \kappa_j(\pi)
}
{
\frac{\theta}{\beta}+ i+1 - \sum_{j=k+1-i}^k \kappa_j(\pi)
}
&
=
\E{
\frac{1}{\frac{\theta}{\beta}+1}
\prod_{i=1}^{k-1}
\frac{
i+1 - \sum_{j=k+1-i}^k \kappa_j(\pi)
}
{
\frac{\theta}{\beta}+ i+1 - \sum_{j=k+1-i}^k \kappa_j(\pi)
}
}
,
\end{align*}
where the expectation is taken relative to $\pi$. By recognizing that the Laplace-Stieltjes transform of $X \sim \mathsf{Exp}(\lambda)$ is $\E{e^{-\theta X}} = \lambda \slash (\theta + \lambda)$, we can further reduce to 
\begin{align*}
\E{
\frac{1}{\frac{\theta}{\beta}+1}
\prod_{i=1}^{k-1}
\frac{
i+1 - \sum_{j=k+1-i}^k \kappa_j(\pi)
}
{
\frac{\theta}{\beta}+ i+1 - \sum_{j=k+1-i}^k \kappa_j(\pi)
}
}
&
=
\E{
e^{-\frac{\theta}{\beta}\sum_{i=1}^{k}T_{\pi,i}}
}
,
\end{align*}
where now the expectation is taken relative to both $\pi$ and the random variables $T_{\pi,i}$. Hence, removing the conditioning on $N = k + 1$, we achieve the stated equivalence.
\hfill \Halmos \endproof

As direct corollaries of this result, we can see how the distribution of objects like the duration, the intensities, and the full arrival epoch sequence depend on the parameters $\alpha$ and $\beta$. For example, by consequence of Theorem~\ref{markovDistThm}, if $0 < \alpha_1 < \beta_1$ and $0 < \alpha_2 < \beta_2$ with $\alpha_1 \slash \beta_1 = \alpha_2 \slash \beta_2$ are parameter pairs defining two UHP models with respective duration random variables $\tau_1$ and $\tau_2$, then $\beta_1 \tau_1$ and $\beta_2 \tau_2$ are equivalent in distribution. Similarly said, for $\rho = \alpha \slash \beta$ held fixed, Theorem~\ref{markovDistThm} shows that the only dependence $\tau$ has on $\beta$ is through the $1 \slash \beta$ coefficient. This realization can actually be traced back to Lemma~\ref{Gentimechange}, implying that this same proportionality can be extended to any arrival epoch, not only the final one. Conditioning yields similar insights even when comparing different values of $\rho$. For example, the conditional LST found in the proof shows that $\beta_1 \tau_1 \mid N_1 = k$ is equivalent to $\beta_1 \tau_2 \mid N_2 = k$ for all $k \in \mathbb{N}$ even if $\alpha_1 \slash \beta_1 \ne \alpha_2 \slash \beta_2$. 

\section{A Modular Simulation Algorithm for General Hawkes Processes}\label{simSec}

In addition to analytic results, the conditional uniformity and parking function decomposition also enable us to propose a novel simulation algorithm for Hawkes process clusters. This new procedure draws upon the combinatorial structure and group theoretic perspectives of  parking functions, mirroring what we have reviewed in Section~\ref{backgroundSec}. As we will see, this new algorithm offers competitive efficiency relative to widely used procedures in the literature. However, like in the preceding demonstration \edit{of analysis}, the most important implication of this methodology again lies in the conditioning. The method we propose first generates the size of the cluster and then generates the arrival epochs accordingly, and we will discuss how this leads to natural applications for rare event simulation and importance sampling.

\begin{algorithm}
\SetAlgoLined

\textbf{Input:} Integrated excitation kernel $G(x) = \int_0^x g(u) \mathrm{d}u$, with $\rho = \int_0^\infty g(u) \mathrm{d}u$

\textbf{Output:} {Cluster of $N$ self-excited arrival epochs $0 = A_0 < A_1 < \dots < A_{N-1}$.}

\begin{enumerate}
\item Generate the total number of events, $N \sim \mathsf{Borel}(\rho)$.
\item Generate the conditionally uniform compensator points, $\Lambda \in \mathbb{R}_+^{N-1}$:
\begin{enumerate}
\item Draw $\pi_i \stackrel{\mathsf{iid}}{\sim} \mathsf{Uni}\{1,N\}$ for $i \in \{1, \dots, N-1\}$, define $\sigma \in \{0,1\}^{N}$ s.t.~$\sigma = [0 \dots 0]$.
\item \textbf{for} $i \in \{1, \dots, N-1\}$:
\item \qquad \textbf{if} $\sigma_{\pi_i} = 0$, set $\sigma_{\pi_i} = 1$; \textbf{else} set $\sigma_{j^*~\mathsf{mod}~N} = 1$, where $j^* = \min\{j > \pi_i \mid \sigma_{j~\mathsf{mod}~N} = 0\}$. 
\item Set $\pi = (\pi - \ell)~\mathsf{mod}~N$ where $\ell$ is s.t.~$\sigma_\ell = 0$.
\item Draw $U_i \stackrel{\mathsf{iid}}{\sim} \mathsf{Uni}(0,1)$  for $i \in \{1, \dots, N-1\}$ and set $\Lambda = \rho\cdot\mathbf{sort}(\pi - U)$.
\end{enumerate}
\item \textbf{return} $A_1, \dots, A_{N-1}$ as the unique solution to the system of equations, $\Lambda_i = \sum_{j=0}^{i-1}G(A_i - A_j)$ for each $i \in \{1, \dots, N-1\}$.
\end{enumerate}
 \caption{Simulation of a Hawkes Process Cluster}
 \label{algSim}
\end{algorithm}

We present the pseudo-code \edit{for} our procedure now in Algorithm~\ref{algSim}. Let us walk through the idea of this method. The algorithm terminates for any Hawkes process with $\rho < 1$, and it uses $G(\cdot)$ and $\rho$ as arguments. In the first step, the size of the cluster $N$ is sampled from the Borel distribution described in Equation~\eqref{borelDist}. All of the remaining steps will then use $N$ as if it was a parameter. Steps 2 and 3 presume that $N > 1$; whenever Step 1 yields $N=1$, the simulation trivially completes with $A_0 = 0$. The second step uses Theorem~\ref{distThm} to generate the vector of compensator transformed arrival epochs, $\Lambda$, by generating a uniformly random parking function, $\pi$. Specifically, Steps 2(a)-(d) employ a group theoretic sampling of $\pi$ that follows directly from Pollak's circle argument from the proof of Lemma~\ref{pollak}. Like what is shown in Figure~\ref{pollakfig}, Step 2(a) generates the preference vector on spaces $1$ through $N$, and Steps 2(b) and (c) ``park'' the preferences on the circle, as recorded in the vector $\sigma$. Then, Step 2(d) identifies the empty space on the circle and rotates it so that the open spot becomes space $N$, yielding $\pi$ as a parking function of length $N-1$. Step 2(e) then implements Theorem~\ref{distThm}, and Step 3 returns the unique vector of arrival epochs associated with this vector of compensator points. This brings us to our third main result.

\begin{proposition}\label{algprop}
The output of Algorithm~\ref{algSim} exactly follows the joint distribution of the times and size of a Hawkes process cluster.
\end{proposition}

The proof of Proposition~\ref{algprop} follows directly from Theorem~\ref{distThm} and the proof of Lemma~\ref{pollak}. One of the benefits of this algorithm is that it is modular, in the sense of the word as it might be used by enthusiasts of high fidelity audio equipment or designers of building construction. That is, it can be built from the pieces one prefers. Each of the three primary steps can be conducted through one's selected method, swapping a style of one component out for another approach as desired. For example, our Step 2 uses Pollak's circular parking argument to generate $\pi$, which aligns with \citet{kenyon2021parking}. However, one could just as readily use a similar idea described in \citet{diaconis2017probabilizing}, where the preference vector is simply iteratively increased by 1 modulo $N$ until it is a valid parking function. Similarly, one could forgo Theorem~\ref{distThm} and instead sample from the polytope $\mathcal{P}_k$ in Step 2 through use of Markov Chain Monte Carlo (MCMC) methods like Metropolis-Hastings or hit-and-run algorithms \citep[see, e.g.,][and references therein]{chen2018fast}. Likewise, our implementations of Algorithm~\ref{algSim} have generated $N$ in Step 1 directly through the probability mass function in~\eqref{borelDist}, but one could, for example, instead simulate a Poisson branching process to yield the same distribution. 

This speaks to another advantage to which we have already alluded: the potential for application in rare event simulation and importance sampling. Algorithm~\ref{algSim} untangles the cluster distribution and the cluster size. That is, this method allows one to specify a particular value or, more generally, a target distribution of the cluster size in Step 1 and then generates a collection of cluster arrival epochs matching this desired cardinality in Steps 2 and 3. We expect this to be of particular practical value, as the Borel distribution might be called \emph{nearly} or \emph{asymptotically} heavy tailed. As $\rho \to 1$, Equation~\eqref{borelDist} yields a valid probability distribution on the positive integers with no finite moments. Hence, it is quite natural for a Hawkes process with $\rho$ near 1 to experience clusters of extreme size, and Algorithm~\ref{algSim} provides a controlled way of simulating and evaluating these scenarios. To the best of our knowledge, no prior Hawkes process simulation algorithm allowed for user selection of the cluster size, as this is a direct consequence of the conditional uniformity we have explored throughout this paper.

\begin{table}[h]
\renewcommand{\arraystretch}{1.25}
\caption{Observed run times (seconds) of Hawkes process cluster simulation algorithms across various choices of excitation kernels ($\boldsymbol{2^{20}}$ replications).}\label{simTable}
\centering
\begin{small}
\centering
\begin{tabular}{c | c c c c c c c c}
& \multicolumn{8}{c}{\textbf{Excitation Kernel} $\boldsymbol{g(x)}$}
\\
\hline
\textbf{Simulation Procedure} 
& $3 e^{-4 x}$ & $15 e^{-16 x}$ & $63 e^{-64 x}$ & $255 e^{-256 x}$
& $\frac{1}{(2+x)^2}$ & $\frac{3}{(4+x)^2}$ & $\frac{7}{(8+x)^2}$ & $\frac{15}{(16+x)^2}$
\\
\hline
Algorithm~\ref{algSim} & 5.6  &  11.4  & 30.6 & 94.4 & 7.4 & 15.2 & 39.1 & 143.6 
\\
\citet{dassios2013exact} & 2.1 & 12.0 & 47.1 & 183.2 & $\times$ & $\times$ & $\times$ & $\times$
\\
\citet{hawkes1974cluster} & 42.2 & 182.0 & 727.8 & 2,945.5 & 19.3 & 42.4 & 88.7 & 181.5 
\\
\end{tabular}
\end{small}
\end{table}

To measure the efficiency and performance of Algorithm~\ref{algSim} relative to the literature, in Table~\ref{simTable} we compare the observed run times of this procedure against the literature.\footnote{\edit{The observed KS statistics between the empirical distributions of $\tau$ and $N$ from Algorithm~\ref{algSim} and those from the literature are no more than 0.001 in the exponential kernel experiment, and not more than 0.006 for the power law.}} Specifically, we consider two prior algorithms: the classical non-stationary Poisson clustering process perspective provided by \citet{hawkes1974cluster} and the exact simulation method by \citet{dassios2013exact}. Let us note that other well known Hawkes process simulation methods, such as the thinning procedure from \citet{lewis1979simulation,ogata1981lewis} or Algorithm 7.4.III from \citet{daley2003introduction}, cannot be applied to \edit{the} generation of clusters alone, as they implicitly rely on the presence of a continual baseline stream. \citet{hawkes1974cluster}'s method is essentially exactly like the corresponding definition: for each arrival, simulate its own descendants according to a non-stationary Poisson process with rate given by the kernel $g(\cdot)$, and repeat this process until there are no more descendants. The generality of this procedure has made it quite popular; it is a backbone subroutine of the perfect Hawkes process simulation algorithms in \citet{moller2005perfect,chen2021perfect}, and~\citet{chen2020perfect}. The \citet{dassios2013exact} procedure (Algorithm 3.1 in that paper) only applies to the exponential kernel case and it relies on that Markov assumption. At each step, a weighted coin is flipped to decide if there will be another arrival given the current value of the intensity. If there will be another, then the time and intensity are updated accordingly, and otherwise the simulation terminates.

We conduct the simulation in Table~\ref{simTable} with two families of excitation kernels, the Markovian exponential decay kernel and the heavy-tailed power law kernel, which is often called ``Omori's law'' in seismology \citep{ogata1998space}. Specifically, we consider a series of kernels in each family with increasing mean cluster size: $g(x) = (4^m - 1)e^{-4^m x}$ for the exponential, which will have mean cluster size $4^m$, and $g(x) = (2^m - 1)\slash(2^m + x)^2$ for the power law, which will have mean cluster size $2^m$. All three algorithms can be used to simulate clusters under the exponential kernel, but only \citet{hawkes1974cluster} can be compared to Algorithm~\ref{algSim} outside of the Markovian setting.  Table~\ref{simTable} shows that the speed of Algorithm~\ref{simTable} is broadly competitive with the literature and superior in many cases. In particular,  in the case of the exponential kernel, the parking function simulation of the Hawkes process proves to be more efficient than either alternative as the clusters grow large. Here, Section~\ref{markovSec} shows how Step 3 can be solved in closed form, and so the efficiency of generating one Borel distributed random variable and one random parking function outpaces the performance of either simulating a non-stationary Poisson process many times in \citet{hawkes1974cluster} or performing many iterative computations of the intensity and the next event epoch in \citet{dassios2013exact}. When the clusters are mostly small though, the set-up of the circular parking function sampling is likely less efficient than \citet{dassios2013exact}'s simple steps, and so Algorithm~\ref{algSim} only outperforms \citet{hawkes1974cluster} in such settings. For general decay kernels, closed form solutions to Step 3 may not be available and one might instead turn to root-finding methods, like what we have done here for the power law kernel. Here we make relatively naive use of Newton's method and, while the efficiency of the Borel and parking function generations helps Algorithm~\ref{algSim}, it is clear that its performance is bottlenecked by the Newton calculations. Thus, as the cluster size increases our method should eventually be outpaced by the \citet{hawkes1974cluster} approach. Hence, careful design of Step 3 of Algorithm~\ref{algSim} is an interesting future direction for this work, and it seems likely that more efficient approaches are available when tailoring to a particular kernel. As we have discussed, though, we anticipate the primary benefit of this algorithm to lie in the conditioning on $N$, since the size distribution's propensity for large values creates an opportunity for rare event simulation.

\section{Discussion and Conclusion}\label{concSec}

In this paper, we have seen that uniformly random parking functions constitute hidden spines \edit{in} Hawkes process clusters, providing a decomposition structure for the conditionally uniform compensator transform of the cluster's arrival epochs. Hence, from the cluster size and a random parking function, we can faithfully reconstruct the full sequence of events in the cluster. We have also demonstrated the impact of this connection \edit{on} both analysis and simulation. For the former, we have found a surprising distributional equivalence between the duration of the Markovian Hawkes cluster and a random sum of \edit{conditionally independent} exponential random variables. For the latter, we have demonstrated that the resulting sampling methodology offers both efficiency and, perhaps more importantly, fine control over the experiment through cluster-size-conditioning.

Much of our discussion of the connection between Hawkes clusters and parking functions has centered around what the discrete object can do for the continuous, but we can also see that there is some partial reciprocity available. That is, it is an immediate consequence of Theorem~\ref{distThm} that, for the arrival epochs drawn from any Hawkes process cluster, the shuffled ceiling of the compensator points will be a parking function, where only the length of the parking function will be dependent on the parameters of the Hawkes process. This underscores a point that we have mentioned only briefly, which is that Lemma~\ref{Gentimechange} shows how one Hawkes process cluster can readily be transformed to another cluster driven by an entirely different excitation kernel so long as the values of $\rho$ match. Hence, one could replace Steps 1 and 2 in Algorithm~\ref{algSim} with simulating the cluster according to some efficient Hawkes process procedure, such as \citet{dassios2013exact}, and then transform to the true targeted excitation function in Step 3. In a similar notion, one could also petition to the classical random time change theorem in place of Lemma~\ref{Gentimechange} and replace Steps 1 and 2 with a unit\edit{-}rate Poisson process. That is, one could run a unit\edit{-}rate Poisson process until the first index $i \geq 1$ such that the Poisson arrival epoch $T_i$ exceeds $i\rho$, and then return times 1 through $i-1$ as the compensator points and proceed to Step 3. Analogously by the closure of the exponential distribution when multiplying by a constant, one could also run a Poisson process at rate $\rho$ and then compare $T_i$ just to the index $i$. This is close to the idea of Algorithm 7.4.III from \citet{daley2003introduction}, but the above proposal would terminate with an absorption event for the end of the cluster. By comparison, it is inherent to Algorithm 7.4.III that the point process continues indefinitely and that the compensator grows without bound, as otherwise ``the final interval is then
infinite and so cannot belong to a unit-rate Poisson process'' \citep{daley2003introduction}. In numerical experiments, we found that the speed of this alternate method is essentially identical to the \citet{dassios2013exact} algorithm. This makes sense, because they are mostly doing the same thing: generating exponentials until an absorption event occurs. Hence, there could be some efficiency gains over our algorithm when $\rho$ is \edit{small enough}, but regardless, like  \citet{dassios2013exact}, this approach does not offer the control over the support and distribution of $N$ like Algorithm~\ref{algSim} does.

In the broader perspective of the Hawkes process literature, we have studied the same model as contemporary works like \citet{chen2021perfect,graham2021regenerative}, which is a linear \citep[relative to the non-linear generalization introduced by][]{bremaud1996stability} univariate \citep[relative to the mutually-exciting version actually dating back to][]{hawkes1971spectra} Hawkes process with general excitation kernel \citep[relative to kernel-specific works like][]{dassios2013exact}. For generalizations of these results, our initial suspicion is that multiple dimensions is the most promising next step, as there are analogs of much of the background fundamentals from which we have drawn. In particular, the random time change does still hold in higher dimensions, where we can build from both results for general processes \citep[e.g.][]{brown1988simple} and for Hawkes processes specifically \citep[e.g.][]{embrechts2011multivariate}. One can also leverage similar time-agnostic multi-type branching process perspectives for the cluster size distribution  \citep[e.g.][]{good1960generalizations}. In fact, a very similar proof to Lemma~\ref{Gentimechange} can create conditionally uniform analogs of the two styles of multivariate random time change in \citet{embrechts2011multivariate}. The remaining challenge arises while converting from conditionally uniform compensator points to the arrival epochs in each sub-stream, as it may be possible that there are multiple solutions to the compensator system of equations. Hence, we are quite interested in this future direction, particularly given the variety of recent interest in these models, such as in works like \citet{ait2015modeling,nickel2020learning,daw2021co,karim2021exact}.

\edit{Another intriguing open problem arises in trying to transport these results from conditioning on the long-run total size of the cluster to instead conditioning on the number of events by some given time $T$. A transient result like this would likely be of broad interest. Arguments similar to that of Lemma~\ref{Gentimechange} should produce a transient result like
\begin{align*}
f(\Lambda_1, \Lambda_2, \dots, \Lambda_k \mid N_T=k+1)
&=
\frac{e^{-\Lambda(T)}}{\PP{N_T = k+1}}
,
\end{align*}
however this expression may be deceptively simple. Setting aside the fact that the distribution of $N_T$ may not be as accessible as that of $N$, what this form of the density hides is the dependence of $\Lambda(T)$ on $\Lambda_1$ through $\Lambda_k$. That is, $\Lambda(T) = \sum_{i=0}^k G(T-A_i)$, and the epochs, $A_1, \dots, A_k$, are deterministic functions of the compensator points, $\Lambda_1, \dots, \Lambda_k$. All this is to say, the conditional density is no longer uniform. While each $\Lambda_i$ will continue to satisfy the constraints $\Lambda_{i-1} < \Lambda_i < i\rho$ (as well as the added constraint $\Lambda_i < \Lambda(T)$), changes to the values of these points will change $\Lambda(T)$ and thus change the density. This suggests a direction that may offer more immediate tractability. Rather than fixing the ending time $T$, one can instead fix the ending compensator value $\Lambda(T)$. In this case, the conditional density should again be uniform, and the compensator points should lie on a polytope that is the intersection of $\mathcal{P}_k$ and the hypercube where every coordinate is between 0 and $\Lambda(T)$. While this result would leverage the transience in the transformed compensator space, it is of course more desirable to analyze true transience in time. Bridging these perspectives stands as a highly interesting direction of future work.}

Finally, to that end, let us emphasize that this paper's analysis has not been a purely theoretical exercise. Our original motivation was to study Hawkes process clusters in an operational context inspired by the data and application in \citet{daw2021co}, but we found the methodological cupboard not yet full enough for the questions we sought to answer. Hence, we are quite interested in returning to these problems with the ability to address the cluster duration and chronology now in hand. Building from the model of \citet{daw2021co}, the Hawkes cluster duration can be seen as the length of a co-produced service exchange, \edit{and, more specifically, the cluster epochs mark the contributions and  points of interactions between the customer and service agent. In the context of \citet{daw2021co}, these epochs are the timestamps of messages exchanged within a text-based contact center's conversational service}. Understanding \edit{these distributions} allows us to evaluate how the history-dependent service process will impact the service system overall, leading to new possibilities in the analysis of natural queueing theoretic questions like staffing and routing decisions. \edit{The service time distribution is a cornerstone of any queueing model, and a myriad, if not virtually all, queueing formulas depend on the mean of this distribution (at the very least). In services modeled as a Hawkes cluster like in \citet{daw2021co}, both this mean and this distribution have been out of reach even for the simplest forms of this process, like \citet{moller2005perfect} described. Thanks to the parking function structure uncovered in Theorems~\ref{distThm} and~\ref{markovDistThm}, this is no longer the case. Leveraging this hidden spine} is of foremost intrigue to us as a direction of future research, and we \edit{anticipate} that the conditional uniformity property and parking function decomposition \edit{can} provide valuable insight into \edit{these classic service-side queueing problems.}

\edit{Of course, the Hawkes process finds many an application in operations research models beyond service durations and customer-agent interactions. In fact, well-known examples of the Hawkes show that clusters are both quantities of operational interest in their own right and key components embedded within any use of the point process. In addition to our original motivating application, this stochastic model can be found among representations of financial contagion and risk, product management, customer arrivals, leadership and communication patterns, digital marketing conversions, and -- lest we forget -- outbreaks in pandemics. By the very nature of the Hawkes process and its self-excited epochs, clusters lie at the core of each of these applications. Let us detail them.} 

\edit{For example, \citet{azizpour2018exploring} uses the Hawkes events to model corporate defaults. Here, a cluster contains the collection of firms that fail as a downstream result of one initial default; likewise, the duration captures how long the contagion lasts from first default to last related casualty. Similarly, \citet{mukherjee2022hiding} models the timing of product recalls, and thus each cluster contains those customers who return their purchases because of the influence of returns by their peers. Elsewhere in the transaction timeline, \citet{xu2014path} uses the Hawkes process as a model of customers' online shopping activity under digital marketing campaigns, with the cluster marking the clicks on the path of conversion from advertisement to purchase. Similarly, on the arrival side of queueing, the Hawkes has served as the arrival stream to both single \citep{chen2021perfect} and infinite server queueing models \citep{gao2018functional,daw2017queues,koops2017infinite}, and in these cases a cluster represents the lineage of customers who choose to patron the service because they saw others do it first. \citet{fox2016HawkesEmails} models the communication (and, within this, leadership) structure in an organization through Hawkes processes, and here clusters arise in the response threads and email chains. Finally, when modeling the spread of COVID-19 and other infectious diseases, like what is done in \citet{bertozzi2020challenges}, the Hawkes clusters naturally capture the trace of contacts who pass the sickness on to one another throughout time.}

\edit{While these and many other Hawkes process applications may feature both self-excited arrivals and exogenously driven baseline events, let us recall that the \citet{hawkes1974cluster} representation shows us how the cluster perspective persists. That is, the baseline stream can be thought of as the arrival process \textit{of clusters}. Hence, the model studied in this paper descends off from each of these initial points, and these initial points themselves form a Poisson process that is independent from the subsequent activity within the clusters. So, the distribution of the cluster epochs and durations constitute the nature of the offsets from the well-understood Poisson stream.}

\edit{From what we have seen, a parking function or Dyck path may not have a direct interpretation for these applications. However, it does provide a structure through which the Hawkes cluster becomes easier to understand, and thus easier to use. Through another level of conditioning, one can apply both traditional Poisson conditional uniformity and the Hawkes variant which we have seen here to analyze, say, customer purchasing patterns or COVID-19 outbreaks and glean insight into the distributions, relationships, and dependencies within. These hidden spines simplify the hallmark self-exciting structure of the Hawkes process, offering  answers to fundamental questions about the model. Returning now to the service operations domains from which this question originally sprung, we look forward to leveraging our new knowledge and exploring what more we can learn once we know how long a cluster will last.}

\section*{Acknowledgements} We are grateful for valuable comments and suggestions given by David Eckman, Sam Gutekunst, Shane Henderson, and Galit Yom-Tov throughout this work. Furthermore, we appreciate the insights and feedback provided from the referees, associate editor, and area editor Jose Blanchet, which have likewise improved this paper. We are also grateful to Ozan Bayiz for conducting auxiliary numerical experiments in related explorations.

\bibliographystyle{informs2014} 
\bibliography{Bibliography.bib}

\begin{thebibliography}{54}
\providecommand{\natexlab}[1]{#1}
\providecommand{\url}[1]{\texttt{#1}}
\providecommand{\urlprefix}{URL }

\bibitem[{A{\"\i}t-Sahalia et~al.(2015)A{\"\i}t-Sahalia, Cacho-Diaz,
  \protect\BIBand{} Laeven}]{ait2015modeling}
A{\"\i}t-Sahalia Y, Cacho-Diaz J, Laeven RJ (2015) Modeling financial contagion
  using mutually exciting jump processes. \emph{Journal of Financial Economics}
  117(3):585--606.

\bibitem[{Azizpour et~al.(2018)Azizpour, Giesecke, \protect\BIBand{}
  Schwenkler}]{azizpour2018exploring}
Azizpour S, Giesecke K, Schwenkler G (2018) Exploring the sources of default
  clustering. \emph{Journal of Financial Economics} 129(1):154--183.

\bibitem[{Bertozzi et~al.(2020)Bertozzi, Franco, Mohler, Short,
  \protect\BIBand{} Sledge}]{bertozzi2020challenges}
Bertozzi AL, Franco E, Mohler G, Short MB, Sledge D (2020) The challenges of
  modeling and forecasting the spread of {COVID-19}. \emph{Proceedings of the
  National Academy of Sciences} 117(29):16732--16738.

\bibitem[{Br{\'e}maud(1975)}]{bremaud1975extension}
Br{\'e}maud P (1975) An extension of {Watanabe's} theorem of characterization
  of {Poisson} processes over the positive real half line. \emph{Journal of
  Applied Probability} 12(2):396--399.

\bibitem[{Br{\'e}maud(1981)}]{bremaud1981point}
Br{\'e}maud P (1981) \emph{Point processes and queues: martingale dynamics},
  volume~50 (Springer).

\bibitem[{Br{\'e}maud \protect\BIBand{}
  Massouli{\'e}(1996)}]{bremaud1996stability}
Br{\'e}maud P, Massouli{\'e} L (1996) Stability of nonlinear {Hawkes}
  processes. \emph{The Annals of Probability} 1563--1588.

\bibitem[{Br{\'e}maud et~al.(2002)Br{\'e}maud, Nappo, \protect\BIBand{}
  Torrisi}]{bremaud2002rate}
Br{\'e}maud P, Nappo G, Torrisi GL (2002) Rate of convergence to equilibrium of
  marked {H}awkes processes. \emph{Journal of Applied Probability}
  39(1):123--136.

\bibitem[{Brown et~al.(2002)Brown, Barbieri, Ventura, Kass, \protect\BIBand{}
  Frank}]{brown2002time}
Brown EN, Barbieri R, Ventura V, Kass RE, Frank LM (2002) The time-rescaling
  theorem and its application to neural spike train data analysis. \emph{Neural
  Computation} 14(2):325--346.

\bibitem[{Brown et~al.(2005)Brown, Gans, Mandelbaum, Sakov, Shen, Zeltyn,
  \protect\BIBand{} Zhao}]{brown2005statistical}
Brown L, Gans N, Mandelbaum A, Sakov A, Shen H, Zeltyn S, Zhao L (2005)
  Statistical analysis of a telephone call center: A queueing-science
  perspective. \emph{Journal of the American Statistical Association}
  100(469):36--50.

\bibitem[{Brown \protect\BIBand{} Nair(1988)}]{brown1988simple}
Brown TC, Nair MG (1988) A simple proof of the multivariate random time change
  theorem for point processes. \emph{Journal of Applied Probability}
  25(1):210--214.

\bibitem[{Carlsson \protect\BIBand{} Mellit(2018)}]{carlsson2018proof}
Carlsson E, Mellit A (2018) A proof of the shuffle conjecture. \emph{Journal of
  the American Mathematical Society} 31(3):661--697.

\bibitem[{Chen(2021)}]{chen2021perfect}
Chen X (2021) Perfect sampling of {H}awkes processes and queues with {H}awkes
  arrivals. \emph{Stochastic Systems} 11(13):264--283.

\bibitem[{Chen \protect\BIBand{} Wang(2020)}]{chen2020perfect}
Chen X, Wang X (2020) Perfect sampling of multivariate {H}awkes processes.
  \emph{2020 Winter Simulation Conference (WSC)}, 469--480 (IEEE).

\bibitem[{Chen et~al.(2018)Chen, Dwivedi, Wainwright, \protect\BIBand{}
  Yu}]{chen2018fast}
Chen Y, Dwivedi R, Wainwright MJ, Yu B (2018) Fast {MCMC} sampling algorithms
  on polytopes. \emph{The Journal of Machine Learning Research}
  19(1):2146--2231.

\bibitem[{Costa et~al.(2020)Costa, Graham, Marsalle, \protect\BIBand{}
  Tran}]{costa2020renewal}
Costa M, Graham C, Marsalle L, Tran VC (2020) Renewal in {Hawkes} processes
  with self-excitation and inhibition. \emph{Advances in Applied Probability}
  52(3):879--915.

\bibitem[{Daley \protect\BIBand{} Vere-Jones(2003)}]{daley2003introduction}
Daley DJ, Vere-Jones D (2003) \emph{An introduction to the theory of point
  processes: Volume I: Elementary theory and methods} (Springer).

\bibitem[{Dassios \protect\BIBand{} Zhao(2013)}]{dassios2013exact}
Dassios A, Zhao H (2013) Exact simulation of {H}awkes process with
  exponentially decaying intensity. \emph{Electronic Communications in
  Probability} 18.

\bibitem[{Daw et~al.(2021)Daw, Castellanos, Yom-Tov, Pender, \protect\BIBand{}
  Gruendlinger}]{daw2021co}
Daw A, Castellanos A, Yom-Tov G, Pender J, Gruendlinger L (2021) The
  co-production of service: modeling service times in contact centers using
  {Hawkes} processes. \emph{Available at SSRN 3817130} .

\bibitem[{Daw \protect\BIBand{} Pender(2018)}]{daw2017queues}
Daw A, Pender J (2018) Queues driven by {H}awkes processes. \emph{Stochastic
  Systems} 8(3):192--229.

\bibitem[{Diaconis \protect\BIBand{} Hicks(2017)}]{diaconis2017probabilizing}
Diaconis P, Hicks A (2017) Probabilizing parking functions. \emph{Advances in
  Applied Mathematics} 89:125--155.

\bibitem[{Ding et~al.(2009)Ding, Giesecke, \protect\BIBand{}
  Tomecek}]{ding2009time}
Ding X, Giesecke K, Tomecek PI (2009) Time-changed birth processes and
  multiname credit derivatives. \emph{Operations Research} 57(4):990--1005.

\bibitem[{Embrechts et~al.(2011)Embrechts, Liniger, \protect\BIBand{}
  Lin}]{embrechts2011multivariate}
Embrechts P, Liniger T, Lin L (2011) Multivariate {Hawkes} processes: an
  application to financial data. \emph{Journal of Applied Probability}
  48(A):367--378.

\bibitem[{Feller(2008)}]{feller2008introduction}
Feller W (2008) \emph{An introduction to probability theory and its
  applications: Volume 2} (John Wiley \& Sons).

\bibitem[{Foata \protect\BIBand{} Riordan(1974)}]{foata1974mappings}
Foata D, Riordan J (1974) Mappings of acyclic and parking functions.
  \emph{Aequationes Mathematicae} 10(1):10--22.

\bibitem[{Fox et~al.(2016)Fox, Short, Schoenberg, Coronges, \protect\BIBand{}
  Bertozzi}]{fox2016HawkesEmails}
Fox EW, Short MB, Schoenberg FP, Coronges KD, Bertozzi AL (2016) Modeling
  e-mail networks and inferring leadership using self-exciting point processes.
  \emph{Journal of the American Statistical Association} 111(514):564--584.

\bibitem[{Gao \protect\BIBand{} Zhu(2018)}]{gao2018functional}
Gao X, Zhu L (2018) Functional central limit theorems for stationary {H}awkes
  processes and application to infinite-server queues. \emph{Queueing Systems}
  90(1-2):161--206.

\bibitem[{Giesecke \protect\BIBand{} Tomecek(2005)}]{giesecke2005dependent}
Giesecke K, Tomecek P (2005) Dependent events and changes of time.
  \emph{Cornell University Technical Report} .

\bibitem[{Good(1960)}]{good1960generalizations}
Good IJ (1960) Generalizations to several variables of {Lagrange's} expansion,
  with applications to stochastic processes. \emph{Mathematical Proceedings of
  the Cambridge Philosophical Society}, volume~56, 367--380 (Cambridge
  University Press).

\bibitem[{Graham(2021)}]{graham2021regenerative}
Graham C (2021) Regenerative properties of the linear {H}awkes process with
  unbounded memory. \emph{The Annals of Applied Probability} 31(6):2844--2863.

\bibitem[{Hawkes(1971)}]{hawkes1971spectra}
Hawkes AG (1971) Spectra of some self-exciting and mutually exciting point
  processes. \emph{Biometrika} 58(1):83--90.

\bibitem[{Hawkes \protect\BIBand{} Oakes(1974)}]{hawkes1974cluster}
Hawkes AG, Oakes D (1974) A cluster process representation of a self-exciting
  process. \emph{Journal of Applied Probability} 11(3):493--503.

\bibitem[{Karim et~al.(2021)Karim, Laeven, \protect\BIBand{}
  Mandjes}]{karim2021exact}
Karim R, Laeven RJ, Mandjes M (2021) Exact and asymptotic analysis of general
  multivariate {Hawkes} processes and induced population processes. \emph{arXiv
  preprint arXiv:2106.03560} .

\bibitem[{Kenyon \protect\BIBand{} Yin(2021)}]{kenyon2021parking}
Kenyon R, Yin M (2021) Parking functions: From combinatorics to probability.
  \emph{arXiv preprint arXiv:2103.17180} .

\bibitem[{Kim \protect\BIBand{} Whitt(2014)}]{kim2014call}
Kim SH, Whitt W (2014) Are call center and hospital arrivals well modeled by
  nonhomogeneous {P}oisson processes? \emph{Manufacturing \& Service Operations
  Management} 16(3):464--480.

\bibitem[{Konheim \protect\BIBand{} Weiss(1966)}]{konheim1966occupancy}
Konheim AG, Weiss B (1966) An occupancy discipline and applications. \emph{SIAM
  Journal on Applied Mathematics} 14(6):1266--1274.

\bibitem[{Koops et~al.(2018)Koops, Saxena, Boxma, \protect\BIBand{}
  Mandjes}]{koops2017infinite}
Koops D, Saxena M, Boxma O, Mandjes M (2018) Infinite-server queues with
  {H}awkes input. \emph{Journal of Applied Probability} 55(3):920--943.

\bibitem[{Lewis \protect\BIBand{} Shedler(1979)}]{lewis1979simulation}
Lewis PW, Shedler GS (1979) Simulation of nonhomogeneous {P}oisson processes by
  thinning. \emph{Naval Research Logistics Quarterly} 26(3):403--413.

\bibitem[{Meyer(1971)}]{meyer1971demonstration}
Meyer PA (1971) D{\'e}monstration simplifi{\'e}e d'un th{\'e}or{\`e}me de
  {Knight}. \emph{S{\'e}minaire de probabilit{\'e}s de Strasbourg} 5:191--195.

\bibitem[{M{\o}ller \protect\BIBand{} Rasmussen(2005)}]{moller2005perfect}
M{\o}ller J, Rasmussen JG (2005) Perfect simulation of {H}awkes processes.
  \emph{Advances in applied probability} 37(3):629--646.

\bibitem[{Mukherjee et~al.(2022)Mukherjee, Ball, Wowak, Natarajan,
  \protect\BIBand{} Miller}]{mukherjee2022hiding}
Mukherjee UK, Ball GP, Wowak KD, Natarajan KV, Miller JW (2022) Hiding in the
  herd: The product recall clustering phenomenon. \emph{Manufacturing \&
  Service Operations Management} 24(1):392--410.

\bibitem[{Nickel \protect\BIBand{} Le(2020)}]{nickel2020learning}
Nickel M, Le M (2020) Learning multivariate {Hawkes} processes at scale.
  \emph{arXiv preprint arXiv:2002.12501} .

\bibitem[{Oakes(1975)}]{oakes1975markovian}
Oakes D (1975) The {Markovian} self-exciting process. \emph{Journal of Applied
  Probability} 12(1):69--77.

\bibitem[{Ogata(1981)}]{ogata1981lewis}
Ogata Y (1981) On {L}ewis' simulation method for point processes. \emph{IEEE
  Transactions on Information Theory} 27(1):23--31.

\bibitem[{Ogata(1988)}]{ogata1988statistical}
Ogata Y (1988) Statistical models for earthquake occurrences and residual
  analysis for point processes. \emph{Journal of the American Statistical
  association} 83(401):9--27.

\bibitem[{Ogata(1998)}]{ogata1998space}
Ogata Y (1998) Space-time point-process models for earthquake occurrences.
  \emph{Annals of the Institute of Statistical Mathematics} 50(2):379--402.

\bibitem[{Ozaki(1979)}]{ozaki1979maximum}
Ozaki T (1979) Maximum likelihood estimation of {H}awkes' self-exciting point
  processes. \emph{Annals of the Institute of Statistical Mathematics}
  31(1):145--155.

\bibitem[{Reynaud-Bouret \protect\BIBand{} Roy(2007)}]{reynaud2007some}
Reynaud-Bouret P, Roy E (2007) Some non asymptotic tail estimates for {H}awkes
  processes. \emph{Bulletin of the Belgian Mathematical Society-Simon Stevin}
  13(5):883--896.

\bibitem[{Riordan(1969)}]{riordan1969ballots}
Riordan J (1969) Ballots and trees. \emph{Journal of Combinatorial Theory}
  6(4):408--411.

\bibitem[{Stanley(1997)}]{stanley1997parking}
Stanley RP (1997) Parking functions and noncrossing partitions. \emph{The
  Electronic Journal of Combinatorics} R20--R20.

\bibitem[{Stanley(1999)}]{stanley1999enumerative}
Stanley RP (1999) \emph{Enumerative Combinatorics}, volume~2 (Cambridge
  University Press).

\bibitem[{Stanley \protect\BIBand{} Pitman(2002)}]{stanley2002polytope}
Stanley RP, Pitman J (2002) A polytope related to empirical distributions,
  plane trees, parking functions, and the associahedron. \emph{Discrete \&
  Computational Geometry} 27(4):603--602.

\bibitem[{Watanabe(1964)}]{watanabe1964discontinuous}
Watanabe S (1964) On discontinuous additive functionals and {L}{\'e}vy measures
  of a {Markov} process. \emph{Japanese journal of mathematics: transactions
  and abstracts}, volume~34, 53--70 (The Mathematical Society of Japan).

\bibitem[{Xu et~al.(2014)Xu, Duan, \protect\BIBand{} Whinston}]{xu2014path}
Xu L, Duan JA, Whinston A (2014) Path to purchase: A mutually exciting point
  process model for online advertising and conversion. \emph{Management
  Science} 60(6):1392--1412.

\bibitem[{Yan(2015)}]{yan2015parking}
Yan CH (2015) Parking functions. \emph{Handbook of enumerative combinatorics},
  859--918 (Chapman and Hall/CRC).

\end{thebibliography}

%
\end{document}